\documentclass[leqno]{siamltex704} % use larger type; default would be 10pt

\usepackage[utf8]{inputenc} % set input encoding (not needed with XeLaTeX)

\usepackage{float}
\usepackage{amsfonts}
\usepackage{amsmath}
\usepackage{mathtools}
\usepackage{hyperref}
\usepackage[hypcap=true,list=true]{subcaption}
\hypersetup{
    unicode=false,          % non-Latin characters in Acrobat’s bookmarks
    pdftitle={Efficient Mesh Optimization Using the Gradient Flow of the Mean Volume},    % title
    pdfkeywords={mesh} {polyhedron} {smoothing} {regularization}, % list of keywords
    colorlinks=true,       % false: boxed links; true: colored links
    pdfborder=0 0 0,
    linkcolor=blue,          % color of internal links (change box color with linkbordercolor)
    citecolor=blue,        % color of links to bibliography
    filecolor=blue,      % color of file links
    urlcolor=blue           % color of external links
}

\usepackage{rotating}

 \newtheorem{rem}[theorem]{Remark}
\usepackage{geometry} % to change the page dimensions
\usepackage{graphicx} % support the \includegraphics command and options
\usepackage{xcolor}
\usepackage{lineno}
\usepackage{booktabs} % for much better looking tables
\usepackage{array} % for better arrays (eg matrices) in maths
\usepackage{paralist} % very flexible & customisable lists (eg. enumerate/itemize, etc.)
\setdefaultleftmargin{1.6em}{2.2em}{1.87em}{1.7em}{1em}{1em}
\usepackage{verbatim} % adds environment for commenting out blocks of text & for better verbatim
\usepackage{fancyhdr} % This should be set AFTER setting up the page geometry
\usepackage{todonotes}

\newcommand{\cg}{\color{gray}}
\newcommand{\cT}{{\cal{T}}}
\newcommand{\cC}{{\cal{C}}}

\newcommand{\R}{\mathbb{R}}
\newcommand{\N}{\mathbb{N}}
\newcommand{\iq}{\operatorname{iq}}

\newcommand{\Hess}{\operatorname{Hess}}

\newcommand{\SO}{\operatorname{SO}}
\newcommand{\id}{\operatorname{id}}
\newcommand{\im}{\operatorname{Im}}

\newcommand{\vol}{\operatorname{vol}}
\newcommand{\area}{\operatorname{area}}
\def\co{\colon\thinspace}
\def\la{\langle}
\def\ra{\rangle}
\newcommand{\tr}{\operatorname{tr}}
\newcommand{\dfn}{\vcentcolon=}

\title{Efficient Mesh Optimization Using\\ the Gradient Flow of the Mean Volume}
\author{Dimitris Vartziotis\footnotemark[2]\ \footnotemark[3]\ \footnotemark[4] \and Benjamin Himpel\footnotemark[2]}
\begin{document}
%\linenumbers
\maketitle
\renewcommand{\thefootnote}{\fnsymbol{footnote}}
\footnotetext[2]{TWT GmbH Science \& Innovation, Department for Mathematical Research \& Services,
  Ernsthaldenstraße 17,
  70565~Stuttgart, Germany}
\footnotetext[3]{NIKI Ltd.\ Digital Engineering, Research Center,
  205~Ethnikis Antistasis Street, 45500~Katsika, Ioannina, Greece}
\footnotetext[4]{Corresponding author. E-mail address: dimitris.vartziotis@nikitec.gr}
\renewcommand{\thefootnote}{\arabic{footnote}}

\begin{abstract}
The signed volume function for polyhedra can be generalized to a mean volume function for volume elements by averaging over the triangulations of the underlying polyhedron. If we consider these up to translation and scaling, the resulting quotient space is diffeomorphic to a sphere. The mean volume function restricted to this sphere is a quality measure for volume elements. We show that, the gradient ascent of this map regularizes the building blocks of hybrid meshes consisting of tetrahedra, hexahedra, prisms, pyramids and octahedra, that is, the optimization process converges to regular polyhedra. We show that the (normalized) gradient flow of the mean volume yields a fast and efficient optimization scheme for the finite element method known as the geometric element transformation method (GETMe). Furthermore, we shed some light on the dynamics of 
this method and the resulting smoothing procedure both theoretically and experimentally.
\end{abstract}



\pagestyle{myheadings}
\thispagestyle{plain}
\markboth{DIMITRIS VARTZIOTIS AND BENJAMIN HIMPEL}{EFFICIENT MESH OPTIMIZATION USING THE MEAN VOLUME}

%\epigraph{Every mathematician has a secret weapon. Mine is Morse theory.}{Raoul Bott}

\section{Introduction}

In the context of the finite element method, mesh quality affects numerical stability as well as solution accuracy of this method \cite{StrangFix2008}. The geometric element transformation method (GETMe) as a smoothing algorithm for tetrahedral meshes is introduced in \cite{VartziotisWipperSchwald2009}, which was generalized to hybrid meshes in \cite{VartziotisPapadrakakis2013,VartziotisWipperGETMeMixed2009,VartziotisWipper2011Hex,VartziotisWipper2011Mixed,VartziotisWipperPapadrakakis2013}. We consider GETMe as a class of smoothing methods based on simple geometric transformations applied iteratively to all elements individually. Several numerical tests have confirmed that the geometric element transformation given in \cite{VartziotisWipperSchwald2009} and its generalizations reliably and efficiently regularizes
the polyhedron types, which are relevant for the finite element method. In this paper, regularization refers to the convergence of the iteratively transformed polyhedra towards a regular polyhedron, and we ask the reader not to confuse this terminology with the methods in numerical analysis that allow you to deal with ill-conditioned problems. GETMe smoothing also significantly reduces errors in solutions of the finite element method and improves solution efficiency for meshes \cite{VartziotisWipperPapadrakakis2013}. The search for a mathematical proof of its qualitative behavior led us to the discovery of a simple and reasonable quality measure for volume elements optimized by polyhedral generalizations of 
the tetrahedral GETMe algorithm in \cite{VartziotisWipperSchwald2009}. Therefore, the purpose of our work is threefold:
\begin{enumerate}
 \item We introduce the mean volume and turn it into a scaling-invariant quality measure.
 \item We consider the gradient flow of the mean volume and discuss its close relationship to the tetrahedral GETMe algorithm presented in \cite{VartziotisWipperSchwald2009}.
 \item We analyze the new GETMe optimizing and untangling algorithms for volume elements induced by the gradient flow of the mean volume both theoretically and experimentally.
\end{enumerate}

\vspace{2ex}

Mesh smoothing methods can be classified \cite{MeiTipper2013_SmoothingPlanarMeshes,Owen1998,Wilson2011_HexahedralSmoothing} as geometry-based \cite{Field1988,VartziotisWipperSchwald2009}, optimi\-zation-based \cite{Escobar20032775,Branets2005, FreitagJonesPlassman1995,FreitagPlassmann2000,Parthasarathy1991,Shivanna2010,LengZhangXu_NovelGeometricFlowDrivenApproach2012,SastryShontz2012_MeshQualityImprovement,BrewerDiachinKnuppLeurentMelander2003}, physics-based \cite{Shimada2000_Bubbles} and combinations thereof \cite{CanannTristanoStaten1998,Freitag1997,ChenTristanoKwok2003}. In order to be more effective, these methods need to be combined with topological modifications \cite{BossenHeckbert1996,FreitagOllivierGooch1997,KlingnerShewchuk2007}. Local or global Optimization-based methods effectively optimize an objective function measuring the quality of elements or the mesh as a whole. They often lend themselves to untangling algorithms \cite{Knupp2001EWC,Li2000,FreitagPlassmann2000,AmentaBernEppstein1999,
Escobar20032775}. Geometry-based methods like the Laplacian \cite{Field1988} and GETMe \cite{VartziotisWipperSchwald2009} smoothings have the advantage of being exceptionally fast, but their effect on the quality of a mesh is heuristic. Therefore, it came as a surprise, that a minor variation of the GETMe method presented in \cite{VartziotisWipperSchwald2009} generalizes to a local optimization-based mesh smoothing and untangling method for hybrid meshes with an objective function induced by a generalization of the volume function.

After a review of the GETMe algorithm \cite{VartziotisWipperSchwald2009} in \S\ref{SecReview} we introduce the mean volume in \S\ref{SecMeanVolume}. Section~\ref{SecGrad} establishes the intimate relationship of the volume and the GETMe algorithm by way of the gradient field. We leverage this relationship in \S\ref{SecHybrid} in order to give a natural generalization of the GETMe smoothing procedure to other volume elements and hybrid meshes. We discuss the dynamics of these optimization and smoothing methods both from a theoretical and experimental point of view in \S\ref{SecDyna}. Due to their technical nature, we postpone the rigorous computations of the singularities of the gradient of the mean volume to \S\ref{SecRegularity}. We finish with a summary in \S\ref{SecOutlook} and an outlook to global optimization-based GETMe methods.

\section{The tetrahedral GETMe algorithm}\label{SecReview}

In \cite{VartziotisWipperSchwald2009} a powerful heuristic method for smoothing tetrahedral meshes was developed and further generalized to other volume elements in \cite{VartziotisPapadrakakis2013,VartziotisWipperGETMeMixed2009,VartziotisWipper2011Hex,VartziotisWipper2011Mixed,VartziotisWipperPapadrakakis2013}. Let us review the key concept from \cite{VartziotisWipperSchwald2009}.

\subsection{Transformation of a tetrahedron}
Let $\tau=(x_1,x_2,x_3,x_4)^{\mathrm{t}} \in \R^{12}$ denote a tetrahedron with the four
pairwise disjoint nodes $x_i\in\mathbb{R}^3$,
$i\in\{1,\dots,4\}$, which is positively oriented. That is, $\det(D(\tau))>0$ with\begin{equation}
  \label{eq:diffmatrix}
   D(\tau):=(x_2-x_1,x_3-x_1,x_4-x_1)
\end{equation}
representing the $(3\times 3)$-matrix of the difference vectors, which span the
tetrahedron $\tau$. Furthermore, let
\begin{equation}\label{eq:normals}
 \begin{split}
n_1 &:= (x_4-x_2)\times (x_3-x_2) \\
n_2 &:= (x_4-x_3)\times (x_1-x_3) \\ 
n_3 &:= (x_2-x_4)\times (x_1-x_4) \\ 
n_4 &:= (x_2-x_1)\times (x_3-x_1) \,.
\end{split}
\end{equation}
denote the inside oriented face normals of $\tau$.
\begin{figure}[ht]
\def\svgwidth{5cm}
\begin{center}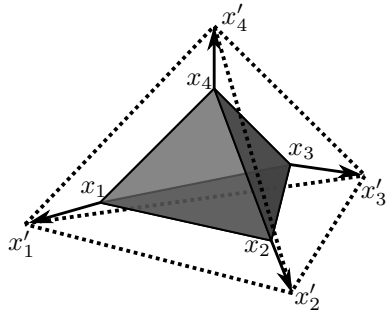\end{center}\caption{Transformation of a tetrahedron\label{Figtransformation}}
\end{figure}

A new tetrahedron $\tau'$ with nodes $x_i'$ is derived from $\tau$
by translating each node $x_i$ using the opposing face normal $n_i$ scaled by
$\sigma/\sqrt{|n_i|}$ for some fixed $\sigma \in \R^+_0$. That is,
\begin{equation}
  \label{eq:transformation}
  \tau'=\begin{pmatrix}
    x_1' \\ x_2' \\ x_3' \\ x_4'
  \end{pmatrix}
  \dfn
\begin{pmatrix}
    x_1 \\ x_2 \\ x_3 \\  x_4 
  \end{pmatrix}+ 
  \sigma
  \begin{pmatrix}
    \frac{1}{\sqrt{|n_1|}}\,n_1 \\
    \frac{1}{\sqrt{|n_2|}}\,n_2 \\
    \frac{1}{\sqrt{|n_3|}}\,n_3 \\
    \frac{1}{\sqrt{|n_4|}}\,n_4
  \end{pmatrix} \,.
\end{equation}

An initial tetrahedron $\tau$ and its transformed counterpart $\tau'$ are depicted in Figure~\ref{Figtransformation}. The edges of the resulting tetrahedron $\tau'$ are indicated by dotted lines.

\subsection{Properties of the transformation}

The quality of tetrahedra is measured by the mean ratio function \cite{K}
\begin{equation}
  \label{eq_meanratio}
  q(\tau):=\frac{3\det (S)^{2/3}}{\|\,S\|_F^2}\,
\end{equation}
with $\|\,S\|_F:=\sqrt{\tr(S^{\mathrm{t}}S)}$ denoting the
Frobenius norm of the matrix $S:=DW^{-1}$. Here $D$ represents the
difference matrix given in \eqref{eq:diffmatrix},
and $W$ denotes the difference matrix of a regular tetrahedron. It holds
that $q(\tau)\in[0,1]$, where very small values indicate nearly degenerate
tetrahedra and larger values almost regular tetrahedra. In particular, it holds
that $q(\tau)=1$, if $\tau$ is regular.
It has been tested heuristically that GETMe increases the mean ratio quality function in \eqref{eq_meanratio}, and that $\sigma$ controls the speed.

\subsection{GETMe smoothing of tetrahedral meshes}\label{SecGETme}

A tetrahedral mesh consists of $n \in \N$ nodes $x_i\in \R^3$ and $m\in \N$ tetrahedra $\tau_j=(x_{i_{j,1}},x_{i_{j,2}},x_{i_{j,3}},x_{i_{j,4}})^{\mathrm{t}}$, where $i_{j,k}\in \{1,\ldots,n\}$, $j=1,\ldots,m$ and $k=1,\ldots,4$. Indices of the nodes for each tetrahedron are ordered such
that all tetrahedra $\tau_j$ are valid in the sense that they are positively
oriented.

Measures for the quality of a tetrahedral mesh are the minimal mean
ratio number
\begin{equation}
  q_{\rm{min}}:=\min_{j\in\{1,\dots, m\}} q(\tau_j)
\end{equation}
and the mean value of all mean ratio numbers
\begin{equation}
  q_{\rm{mean}}:=\frac{1}{m} \sum_{j = 1}^{m} q(\tau_j)\,.
\end{equation}

GETMe smoothing algorithms based on Equation~\eqref{eq:transformation} combine the transformation with a simple scaling and relaxation procedure in order to improve the overall mean ratio mesh quality, while preserving the validity of all tetrahedra after every
iteration step. Scaling is an important aspect of the smoothing procedure, which will leave one of the properties volume, surface area and edge length sum unchanged. There are two flavors:
\begin{itemize}
 \item The simultaneous smoothing applies the transformation simultaneously to all tetrahedra.
 \item The sequential smoothing applies the transformation iteratively to the worst tetrahedron.
\end{itemize}
Usually, the simultaneous and sequential methods have been combined by first simultaneously smoothing the mesh, and then optimizing particularly bad elements by sequential smoothing.

\subsection{Key differences to the present work}

In the following chapters we discuss another geometric element transformation. It differs from Equation~\eqref{eq:transformation} only in its scaling factor, but its GETMe smoothing as described in \S\ref{SecGETme} has the same advantageous behavior. This small change allows for interpreting the transformation as a gradient flow of the volume function (see \S\ref{SecGrad}) and for generalizing it to prisms, pyramids and hexahedra in a natural way (see \S\ref{SecHybrid}). Furthermore, the new point of view shows that our geometric element transformation untangles the individual volume elements (see \S\ref{SecDyna}) and regularizes them (see \S\ref{SecRegularity}).

\section{The mean volume for volume elements}\label{SecMeanVolume}

At first sight, the volume of a polyhedron seems like a rather useless quality measure for polyhedra, because it lacks scaling-invariance and cannot capture information about the (scaling-invariant) shape of the object. It is therefore interesting, that there should be a simple way to turn it into a useful element quality measure for hybrid meshes. First of all, we need to extend the volume to a more abstract notion of volume elements, the building blocks of hybrid meshes. These volume elements are essentially polyhedra, for which we allow non-planar faces and self-intersections. We call this generalization mean volume, and it is defined via triangulations of the underlying polyhedra. This step is followed by a suitable normalization procedure in order to get a scaling-invariant quality measure from the volume.

\subsection{Triangulations}

\begin{figure}[ht]
\def\svgwidth{7cm}
\centering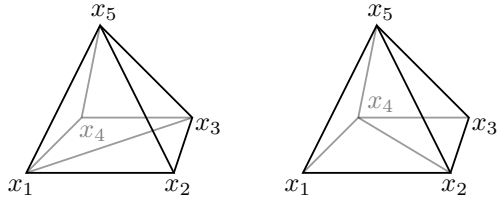\caption{The different triangulations of a pyramid\label{pyramid_triang2}}
\end{figure}

Every convex polyhedron can be triangulated. For example, a pyramid $P=(x_1,\ldots,x_5)$ with apex $x_5$ allows for the two different triangulations depicted in Figure~\ref{pyramid_triang2}.

\begin{definition}\label{DefnTriang}
A triangulation of a convex polyhedron $P$ is given by a family $T=\{\tau_i\}_{i=1,\ldots,m}$, $m \in \N$, of positively oriented tetrahedra $\tau_i$, such that
\begin{enumerate}
 \item $P=\bigcup_i \tau_i$, and
 \item $\tau_i \cap \tau_j$ is either the empty set, a node or a facet for $i\neq j$.
\end{enumerate}
\end{definition}
Observe, that the number $m$ of tetrahedra can be different for different triangulations of the same polyhedron. For example, there are triangulations of a hexahedron with 5 and with 6 tetrahedra.

\begin{figure}[ht]	
\def\svgwidth{14cm}
\begin{center}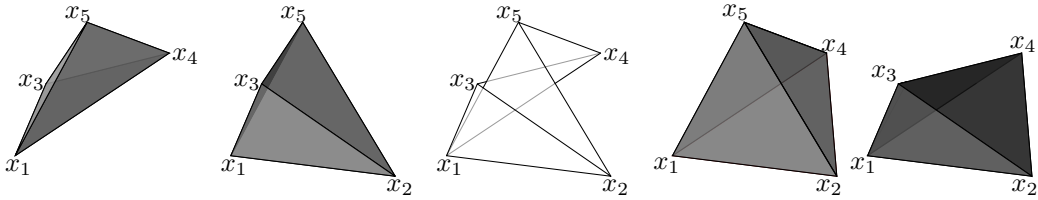\end{center}\caption{A volume element based on the pyramid and its two triangulations (left and right) consisting of two tetrahedra each\label{vol}}
\end{figure}

Volume elements are always associated to some underlying reference polyhedron. This allows us to triangulate any volume element by using the triangulations of the underlying polyhedron. If $x_1,\ldots,x_k$ are the nodes of the volume element, it will be convenient to denote the volume element itself by $x=(x_1,\ldots,x_k)^t \in \R^{3k}$, especially since we will treat $x$ as point on a manifold. 
For tetrahedra, the formal distinction between $\tau$ and $x$ is supposed to clarify, whether the tetrahedron is considered as a volume element or as a polyhedron. While this might seem unnecessary for tetrahedra, it is helpful in the discussion of other volume elements and polyhedra. For example, we have two triangulations of the volume element $x = (x_1,\ldots,x_5)^t$ shown in Figure~\ref{vol}, where the reference polyhedron is the pyramid shown with its triangulations in Figure~\ref{pyramid_triang2}. Notice that only the third tetrahedron in Figure~\ref{vol} is positively oriented, the others are negatively oriented. In particular, there is no canonical three-dimensional geometric object corresponding to this volume element.

Given a volume element $x$, let us denote by $\cT_x$ the set of all triangulations of $x$ corresponding to the triangulations of its underlying polyhedron as described in Definition~\ref{DefnTriang}.

\subsection{The mean volume}

The signed volume of a tetrahedron with vertex coordinates $\tau=(x_1,\ldots,x_4) \in \R^{3\times 4}$ agrees with $1/6$ of the determinant of the difference matrix \eqref{eq:diffmatrix} and can also be written as
\begin{equation}\label{eq_vol}
\vol(\tau) = \frac{1}{6}((x_2-x_1)\times (x_3-x_1))\cdot (x_4-x_1)\;.
\end{equation}
The orientation of the tetrahedron and therefore the sign of the volume function is determined by the order of vertices. Notice that it is therefore not only a well-defined function for valid tetrahedra, but also for invalid and degenerate tetrahedra.

We also have a signed volume function for polyhedra. Given a convex polyhedron $P$, we may consider a triangulation $T$ of $P$ and compute
\begin{equation}\label{eq_volP}
 \vol(P) = \sum_{\tau \in T} \vol(\tau).
\end{equation}

We can also compute this volume function for a volume element. Since the facets are not necessarily planar, this function depends on the triangulation $T$. Given a volume element $x$, we can however simply take the average over all triangulation $\cT_x$.
\begin{definition}\label{DefnMeanVol}
 The mean volume of a volume element $x$ is given by\[
 \vol(x) = \frac{1}{|\cT_x|}\sum_{T\in \cT_x}\sum_{\tau \in T} \vol(\tau).
 \]
\end{definition}
There are alternative definitions for this mean volume function, but this is the most convenient one for our purpose. Notice that there are non-planar polyhedra, for which the mean volume vanishes. For example, the volume element in Figure~\ref{vol} is such an example: only the volume of the third polyhedron is positive, the other three smaller tetrahedra have negative volume.

The mean volume for volume elements generalizes the (generalized) volume for triangulated polyhedral surfaces introduced by Connelly \cite{C}, which was used in the proof of the bellows conjecture \cite{CSW,S}. The quadrilateral (non-planar) facets of the volume elements can be geometrically realized by certain doubly-ruled surfaces \cite{DS}. In applications, we usually only encounter meshes with triangular and quadrilateral facets, and the use of the doubly-ruled surfaces allows us to treat the volume elements as 3--dimensional objects whose volume equals the mean volume. 

\subsection{The quality measure}

Clearly, the mean volume function is translation-invariant, but not scaling-invariant. However, there are several ways to make it scaling-invariant by normalizing it. There is the Frobenius norm used in the mean ratio function \eqref{eq_meanratio}, or we could use a combination of edge length and area. Instead, we define for a volume element $x=(x_1,\ldots,x_k)^t$, for which not all coordinates are equal,
\begin{equation}\label{eq:proj}
\begin{split}
\pi\co \quad\R^{3k} &\to \R^{3k}\\
(x_1,\ldots,x_k)^t &\mapsto \frac{(x_1-x_*,\ldots,x_k-x_*)^t}{\|(x_1-x_*,\ldots,x_k-x_*)^t\|}, \qquad \text{where } x_*=\sum_{i=1}^k x_i.
\end{split}
\end{equation}
The function $\pi$ translates a volume element to its centroid and rescales it. Since $e$ is not entirely degenerate, $(x_1-x_*,\ldots,x_k-x_*)$ does not vanish and $\pi$ is well-defined. The function $q\co \R^{3k}\to \R$ given by \begin{equation}\label{eq:quality}q = \vol \circ \pi\end{equation} is invariant under scaling, translation and rotation, and we will consider it as quality measure for volume elements.

As we will see, regular tetrahedra, hexahedra, octahedra as well as certain symmetric pyramids and prisms maximize this function. Numerical tests and some theoretical evidence suggest, that the global maximum is the only local maximum for all of these volume elements. In particular, optimizing a hybrid mesh with respect to this quality measure should generally improve meshes for its use in finite-element analysis.

\section{The gradient of the mean volume}\label{SecGrad}

Following the gradient flow lines of the mean volume quality measure will certainly improve the volume elements with respect to this quality measure. In this section we will see how the mean volume from \S\ref{SecMeanVolume} relates to the tetrahedral GETMe algorithm described in \S\ref{SecReview}. We find it instructive to introduce a mathematical model for the space of volume elements with a fixed underlying polyhedron, which incorporates the desired invariance under scaling and translation of quality functions.

\subsection{The manifold of volume elements}\label{SecManifold}

Quality measures are often invariant under scaling, translation and rotation. Therefore, we want to find a model for the space of volume elements, on which the analysis of such quality measures is simplified. If we leave out rotations, we will see that the resulting model is simply a sphere.

The space of volume elements $x=(x_1,\ldots,x_k)^t$ is isomorphic to the Euclidean space $\R^{3k}$. Considering these elements up to translation and scaling corresponds to introducing an equivalence relation on $\R^{3k}$ defined by\[
 x \sim x' \vcentcolon\Longleftrightarrow x=\lambda x'+(x_0,\ldots,x_0)^t \text{ for some } \lambda'\in \R\backslash \{0\} \text{ and } x_0\in \R^3\, .
\]
The quotient space by this equivalence relation is not a manifold. If we ignore the degenerate volume elements of the form $(x_0,\ldots,x_0)^t$ we get an open subset of $\R^{3k}$\[
 M \dfn \R^{3k}\backslash\{(x_0,\ldots,x_0)^t\mid x_0\in \R^3\}
\]
whose quotient space $M/\sim$ is simply (diffeomorphic to) a $(3k-4)$--sphere. To be more concrete, 
the map on $M/\sim$ induced by the projection map $\pi$ in \eqref{eq:proj} yields an identification \begin{equation}\label{eq:diffeo}\pi: M/\sim \longrightarrow N\dfn\{x \in \R^{3k} \mid \|x\|=1 \text{ and } x_*\dfn   \sum_{i=1}^n x_i = 0\}.\end{equation}

The differential structure on $N$ induces a differential structure on $M$, and we can consider their tangent bundles $TM$ and $TN$. These are related via the differential $D(\pi)\co TM \to TN$ of $\pi$.

\subsection{The gradient field}\label{sec:gradvol_tet}

Let $C^\infty(M,TM)$ be the space of smooth sections of $TM$, also known as vector fields. For a tetrahedral volume element $x = (x_1,\ldots,x_4)^t \in M$ it is straight-forward to compute from Equation \eqref{eq_vol} that the gradient $\nabla\vol_{\tau}$ of the volume function at $x$ is given by six times the normals $n_i$ in Equation \eqref{eq:normals}, that is,\[
\nabla \vol_{\tau} = \frac{1}{6}\begin{pmatrix}
    n_1 \\
    n_2 \\
    n_3 \\
    n_4
  \end{pmatrix}\,.
\]
The gradient field $\nabla \vol \in C\infty(M,TM)$ of the volume is the vector field given by all gradient vectors $\nabla \vol_x \in T_xM$ of the volume on $M$.

The factors $1/\sqrt{\|n_i\|}$ had been introduced in Equation \eqref{eq:transformation} to ensure scaling invariance. Instead, we use the factor $1/\sqrt{\|n\|}$ for $n = (n_1,\ldots,n_4)^t\in \R^{12}$, which not only ensures scaling invariance, but also preserves the centroid. That is,
\begin{equation}
  \label{eq:transformation2}
  \tau'=\begin{pmatrix}
    x_1' \\ x_2' \\ x_3' \\ x_4'
  \end{pmatrix}
  :=
\begin{pmatrix}
    x_1 \\ x_2 \\ x_3 \\  x_4 
  \end{pmatrix}+ \frac{1}{\sqrt{|n|}}
\begin{pmatrix}
    n_1 \\
    n_2 \\
    n_3 \\
    n_4
\end{pmatrix} \;.
\end{equation}

The intimate relationship between the tetrahedral GETMe transformation and the gradient of the volume function allows us to generalize the tetrahedral GETMe transformation to volume elements such as hexahedra, prisms, pyramids, and even octahedra. This way, we can analyze it not only numerically, but theoretically. In particular, we show that the mean volume quality measure is a height function on the sphere $N$ and discuss some of its properties. Notice that a variation of the necessary scaling mentioned in \S\ref{SecGETme} is even built into the map $\pi$ given in \eqref{eq:proj}. More precisely, we could simply leave $\|(x_1-x_*,\ldots,x_n-x_*)\|$ unchanged.

\section{GETMe smoothing for hybrid meshes}\label{SecHybrid}

In order to give a description of the GETMe smoothing, it will be helpful to discuss the theory in more detail. As we have seen, face normals play an essential role in the tetrahedral GETMe approach. After discussing face normals in more generality, we will explicitly compute gradient vector fields of the mean volume function given in Definition \ref{DefnMeanVol} for pyramids, prisms and octahedra directly from the definition. After some theory about scaling invariance and gradient vector fields, we show how to use the symmetry of Platonic solids to compute the gradient fields for other polyhedra. More precisely, we can use a symmetry-invariant subset of triangulations in order to compute the gradient. For example, the hexahedron has a symmetry-invariant subset of only two triangulations.

\subsection{Face normals}\label{SecNormals}

\begin{figure}[ht]
\def\svgwidth{7cm}
\begin{center}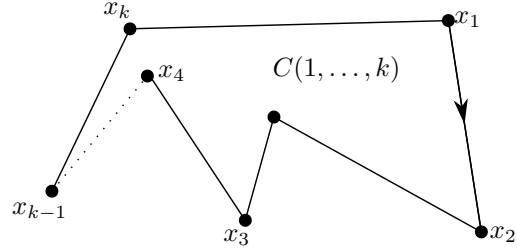\end{center}\caption{The chain $C(1,\ldots,k)$\label{curve}}
\end{figure}

By the face normal of an oriented triangle we simply mean the cross product of two edge vectors compatible with the orientation. Even though cross products are usually applied to edge vectors, we can also apply the cross products to the vertex coordinates themselves, because they are points in the vector space $\R^3$ and therefore vectors themselves. In order to simplify notation, let us define face normals for arbitrary oriented polygonal curves $C(1,\ldots,k)$ in $\R^3$ as shown in Figure~\ref{curve} by
\begin{equation}\begin{split}
\nu(1, \ldots, k)  \coloneqq &(x_2-x_1)\times(x_3-x_1)+(x_3-x_1)\times(x_4-x_1)+\ldots \\
& +(x_{k-2}-x_1)\times(x_{k-1}-x_1) + (x_{k-1}-x_1)\times(x_k-x_1)\\
= &x_{1} \times x_{2} + x_{2} \times x_{3} + \ldots + x_{{k-1}} \times x_{k} + x_{k} \times x_1\;.
\end{split}\end{equation}

\begin{rem}\label{RemIdentities} There are several useful facts about cross products, which can be proven by elementary calculations. We have the identities\begin{align}
(x_{1} -x_{2}) \times (x_{2}-x_{3}) & = \nu(1,2,3),\\
(x_{1}-x_{3}) \times (x_{2}-x_{4}) & = \nu(1,2,3,4),\\
\nu(2,\ldots, k,1)   = \nu(1,\ldots, k)  &  = \sum_{j=2}^{k-1} \nu(1,j,{j+1}). \label{HandyLaw}
\end{align}
The direction of the normal vector $\nu(1,2,3)$ can be found by a right-hand (grip) rule.
Equation~\eqref{HandyLaw} encodes the following properties of cross products:
\begin{enumerate}
\item The face normals of tetrahedra add up to zero.
\item If $C(1,\ldots,k)$ is the boundary of a triangulated oriented surface, then the sum of their normals is independent of the choice of triangulation.
\item If $C(1,\ldots,k)$ is planar and convex, then $\|\nu(1,\ldots,k)\|$ equals twice the area enclosed by $C(1,\ldots,k)$.
\end{enumerate}
\end{rem}

The above remark gives a convenient way of organizing linear combinations of face normals. In particular, it allows us to 
write the gradient field $\nabla \vol \in C^\infty(M,TM)$ of the volume as
\begin{equation}\label{EqTetraX}
\nabla \vol_{\tau} = \frac{1}{6}\begin{pmatrix}
\nu(4,3,2)\\
  \nu(4,1,3)\\
\nu(4,2,1)\\
\nu(1,2,3)
\end{pmatrix}\quad
\text{ for the tetrahedron }\tau=\text{\raisebox{-1cm}{
\def\svgwidth{3cm}
%% Creator: Inkscape inkscape 0.48.3.1, www.inkscape.org
%% Accompanies image file 'tetrahedron.pdf' (pdf, eps, ps)
\begingroup%
  \makeatletter%
  \providecommand\color[2][]{%
    \errmessage{(Inkscape) Color is used for the text in Inkscape, but the package 'color.sty' is not loaded}%
    \renewcommand\color[2][]{}%
  }%
  \providecommand\transparent[1]{%
    \errmessage{(Inkscape) Transparency is used (non-zero) for the text in Inkscape, but the package 'transparent.sty' is not loaded}%
    \renewcommand\transparent[1]{}%
  }%
  \providecommand\rotatebox[2]{#2}%
  \ifx\svgwidth\undefined%
    \setlength{\unitlength}{172.76104736bp}%
    \ifx\svgscale\undefined%
      \relax%
    \else%
      \setlength{\unitlength}{\unitlength * \real{\svgscale}}%
    \fi%
  \else%
    \setlength{\unitlength}{\svgwidth}%
  \fi%
  \global\let\svgwidth\undefined%
  \global\let\svgscale\undefined%
  \makeatother%
  \begin{picture}(1,0.70672248)%
    \put(0,0){\includegraphics[width=\unitlength]{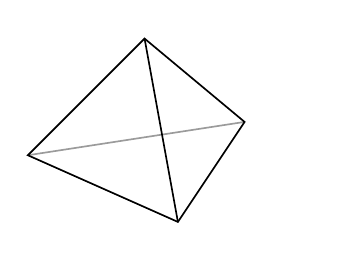}}%
    \put(-0.00845641,0.22909907){\color[rgb]{0,0,0}\makebox(0,0)[lb]{\smash{$x_1$}}}%
    \put(0.45130321,0.02402614){\color[rgb]{0,0,0}\makebox(0,0)[lb]{\smash{$x_2$}}}%
    \put(0.70268273,0.34817328){\color[rgb]{0,0,0}\makebox(0,0)[lb]{\smash{$x_3$}}}%
    \put(0.34545936,0.6293215){\color[rgb]{0,0,0}\makebox(0,0)[lb]{\smash{$x_4$}}}%
  \end{picture}%
\endgroup%
}.}
\end{equation}

It is straight-forward to also write the gradient of the mean volume for other volume elements given in Definition \ref{DefnMeanVol} in terms of face normals. This vector field will yield a flow on the manifold of volume elements. We will see in \S\ref{SecRegularity} and \S\ref{SecDyna} to which extent the optimizing behavior for tetrahedra can be observed and proven for more general polyhedra.

Let us compute the gradient $\nabla\vol$ of the mean volume for the volume element $x$ based on the pyramid. Its two triangulations in $\cT_x$ consist of only two tetrahedra each shown in Figure~\ref{pyramid_triang2}. We compute
\begin{equation}\label{EqGradPyra}
\nabla \vol_x = 
\frac{1}{12}\begin{pmatrix}
\nu(5,4,2) + \nu(5,4,3,2)\\
\nu(5,1,3) + \nu(5,1,4,3)\\
\nu(5,2,4) + \nu(5,2,1,4)\\
\nu(5,3,1) + \nu(5,3,2,1)\\
2 \cdot \nu(1,2,3,4)\end{pmatrix}\;.
\end{equation}

It is also straight-forward to construct the vector field for the prism, because its set of all triangulations has only six elements:
\begin{equation}\label{EqGradPrism}
\nabla \vol_x = \displaystyle \frac{1}{18}\begin{pmatrix}
\nu(3,2,4) + \nu(2,5,4,6,3)\\
\nu(1,3,5) + \nu(3,6,5,4,1)\\
\nu(2,1,6) + \nu(1,4,6,5,2)\\
\nu(5,6,1) + \nu(6,3,1,2,5)\\
\nu(6,4,2) + \nu(4,1,2,3,6)\\
\nu(4,5,3) + \nu(5,2,3,1,4)\\
\end{pmatrix}\quad
\text{ for the prism \raisebox{-1.2cm}{
\def\svgwidth{3cm}
%% Creator: Inkscape inkscape 0.48.3.1, www.inkscape.org
%% Accompanies image file 'prism.pdf' (pdf, eps, ps)
\begingroup%
  \makeatletter%
  \providecommand\color[2][]{%
    \errmessage{(Inkscape) Color is used for the text in Inkscape, but the package 'color.sty' is not loaded}%
    \renewcommand\color[2][]{}%
  }%
  \providecommand\transparent[1]{%
    \errmessage{(Inkscape) Transparency is used (non-zero) for the text in Inkscape, but the package 'transparent.sty' is not loaded}%
    \renewcommand\transparent[1]{}%
  }%
  \providecommand\rotatebox[2]{#2}%
  \ifx\svgwidth\undefined%
    \setlength{\unitlength}{129.63053571bp}%
    \ifx\svgscale\undefined%
      \relax%
    \else%
      \setlength{\unitlength}{\unitlength * \real{\svgscale}}%
    \fi%
  \else%
    \setlength{\unitlength}{\svgwidth}%
  \fi%
  \global\let\svgwidth\undefined%
  \global\let\svgscale\undefined%
  \makeatother%
  \begin{picture}(1,0.87370968)%
    \put(0,0){\includegraphics[width=\unitlength]{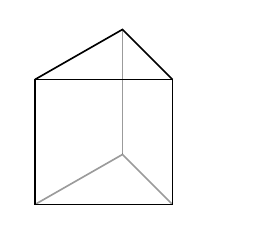}}%
    \put(0.00227698,0.02619827){\color[rgb]{0,0,0}\makebox(0,0)[lb]{\smash{$x_1$}}}%
    \put(0.60812649,0.04256276){\color[rgb]{0,0,0}\makebox(0,0)[lb]{\smash{$x_2$}}}%
    \put(0.31759752,0.32337662){\color[rgb]{0,0,0}\makebox(0,0)[lb]{\smash{$\cg x_3$}}}%
    \put(-0.00922092,0.57948384){\color[rgb]{0,0,0}\makebox(0,0)[lb]{\smash{$x_4$}}}%
    \put(0.67580324,0.56714069){\color[rgb]{0,0,0}\makebox(0,0)[lb]{\smash{$x_5$}}}%
    \put(0.34871962,0.78931113){\color[rgb]{0,0,0}\makebox(0,0)[lb]{\smash{$x_6$}}}%
  \end{picture}%
\endgroup%
}.}
\end{equation}

It is simpler to compute the vector field for the octahedron:
\begin{equation}\label{EqGradOcta}
\nabla\vol_x = \frac{1}{6}\begin{pmatrix}
\nu(2,3,4,5)\\
\nu(1,5,6,3)\\
\nu(1,2,6,4)\\
\nu(1,3,6,5)\\
\nu(1,4,6,2)\\
\nu(2,5,4,3)
\end{pmatrix}
 \quad \text{for \quad
\raisebox{-1cm}{\def\svgwidth{3cm}
%% Creator: Inkscape inkscape 0.48.3.1, www.inkscape.org
%% Accompanies image file 'octahedron.pdf' (pdf, eps, ps)
\begingroup%
  \makeatletter%
  \providecommand\color[2][]{%
    \errmessage{(Inkscape) Color is used for the text in Inkscape, but the package 'color.sty' is not loaded}%
    \renewcommand\color[2][]{}%
  }%
  \providecommand\transparent[1]{%
    \errmessage{(Inkscape) Transparency is used (non-zero) for the text in Inkscape, but the package 'transparent.sty' is not loaded}%
    \renewcommand\transparent[1]{}%
  }%
  \providecommand\rotatebox[2]{#2}%
  \ifx\svgwidth\undefined%
    \setlength{\unitlength}{221.95269775bp}%
    \ifx\svgscale\undefined%
      \relax%
    \else%
      \setlength{\unitlength}{\unitlength * \real{\svgscale}}%
    \fi%
  \else%
    \setlength{\unitlength}{\svgwidth}%
  \fi%
  \global\let\svgwidth\undefined%
  \global\let\svgscale\undefined%
  \makeatother%
  \begin{picture}(1,0.77213982)%
    \put(0,0){\includegraphics[width=\unitlength]{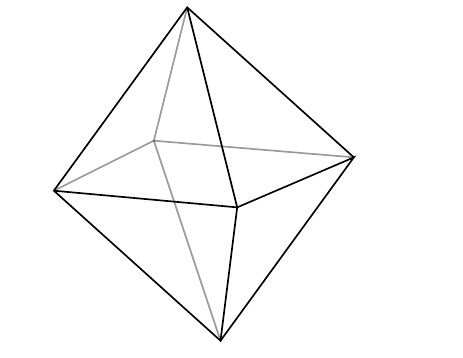}}%
    \put(-0.0065822,0.34188655){\color[rgb]{0,0,0}\makebox(0,0)[lb]{\smash{$x_2$}}}%
    \put(0.52027836,0.2666344){\color[rgb]{0,0,0}\makebox(0,0)[lb]{\smash{$x_3$}}}%
    \put(0.76857752,0.37476555){\color[rgb]{0,0,0}\makebox(0,0)[lb]{\smash{$x_4$}}}%
    \put(0.24866974,0.71189331){\color[rgb]{0,0,0}\makebox(0,0)[lb]{\smash{$x_1$}}}%
    \put(0.31347523,0.0187012){\color[rgb]{0,0,0}\makebox(0,0)[lb]{\smash{$x_6$}}}%
    \put(0.35098267,0.49418555){\color[rgb]{0,0,0}\makebox(0,0)[lb]{\smash{$\cg x_5$}}}%
  \end{picture}%
\endgroup%
}.}
\end{equation}

In order to compute the gradient field for more complicated polyhedra, we will discuss symmetry-invariant subsets of triangulations in \S\ref{SecSymmetry}.

\subsection{Scaling-invariance}

In order for the transformation determined by $\nabla\vol$ to be invariant under scaling, it can be normalized as in \cite{VartziotisWipperSchwald2009} by multiplying each face normal with the inverse of the square root of its norm. Instead, we will generalize the normalization used in Equation \eqref{eq:transformation2}, which preserves the direction of vectors on $M$ and the singularities of the vector field, which makes dynamical studies simpler. 

\begin{definition}\label{DefnNorX} For a vector field $X\in TM$ on $M$ let \[\Psi(X_x)\coloneqq\begin{cases}
\frac{1}{\sqrt{\|X_x\|}}X_x & X_x \neq 0\\
0 & X_x = 0.
\end{cases}
\]
Furthermore, we let $\Psi(X)$ be the vector field given by $\Psi(X)_x = \Psi(X_x)$.
\end{definition}

Scaling and translation invariance of $\Psi(\nabla \vol)$ is equivalent to $\Psi(\nabla \vol)$ being a well-defined vector field on the quotient $M/\sim$.

\begin{lemma}\label{LemYX} The vector field $X=\Psi(\nabla \vol)$
is a vector field on the quotient $M/\sim$.
\end{lemma}

{\em Proof}. Recall the diffeomorphism $\pi\co M/\sim \to N$ defined in \eqref{eq:diffeo}. It suffices to check that $D(\pi)(X_x) = D(\pi)(X_{x'})$ for $x\sim x'$. Notice that $t X_x = X_{t x}$ and $X_x = X_{x + (x_0,\ldots,x_0)}$.   Clearly, if $x' = x + (x_0,\ldots,x_0)$, then
$D(\pi)(X_{x'}) = D(\pi)(X_x)$. If $x' = \lambda x$ with $x_* = \sum x_i = 0$, then using the fact the $\nabla \vol$ preserves the centroid we get \begin{equation*}
D(\pi)(X_{x'}) = \frac{d}{dt} \pi (\lambda x + tX_{\lambda x})|_{t=0} = \frac{d}{dt} \pi (\lambda (x + tX_{x}))|_{t=0} = \frac{d}{dt} \pi (x + tX_{x})|_{t=0} = D(\pi)(X_x).\quad \endproof\end{equation*}
 
\subsection{Gradient fields}\label{SecVecField}

The gradient flow is an important concept and a powerful tool with many applications in topology, global analysis, dynamical systems and mathematical physics. Perelman observed in his proof of Thurston's geometrization conjecture \cite{KL}, that one can interpret Hamilton's Ricci flow as a gradient flow. Morse theory \cite{M,B} with all its generalizations and infinite-dimensional manifestations in mathematical physics like Chern-Simons theory \cite{W}, Yang-Mills theory \cite{AB,A,R} and other gauge theories deduces topological information from the study of the singularities and the flow of a gradient field. It is a natural consequence of the strong results surrounding gradient fields and their flows, that we should make use of them to shed some light on the dynamics on the transformation given in \cite{VartziotisWipperSchwald2009} and generalize it. In this section we will discuss some important properties of gradient fields and Lyapunov functions. Most importantly, we will see, that the (normalized) 
gradient 
flow lines start and 
end at singularities of the gradient field.

Recall the quality measure $q = \vol \circ \pi \co M \to \R$ defined in \eqref{eq:quality}. For $X_x \in T_xN$ we compute for the differential $D$ of the mean volume $D\vol(X_x) = D(\vol|_N)(X_x)$. Then, the gradient field $\nabla (\vol|_N)$ of $\vol|_N$ is implicitly defined using the submanifold metric $\la \cdot,\cdot \ra_N$ on $N\subset M \subset \R^{3k}$ by\[
 D \vol (X_x) = \la \nabla (\vol|_N)_x, X_x \ra_N \quad \text{for all } X_x \in T_xN.
\]
\begin{definition}
 A smooth function $f\co N\to \R$ is a Lyapunov function for a vector field $X \in C^\infty(N,TN)$, if
 \begin{enumerate}
  \item $D f_x (X) > 0$ for all non-singular points $x$,
  \item $D f_x (X) = 0$ if and only if $x$ is a singular point of $X$.
 \end{enumerate}
\end{definition}
A special case of a Lyapunov function is a smooth map $f\co N \to \R$ and its gradient field $\nabla f$, because $Df_x(\nabla f) = \| \nabla f_x\|^2_N >0$ for all non-singular points $x$ and $Df_x(X) =\la\nabla f_x,X\ra_N = 0$ for all $X \in T_x N$.
Since $N$ is a compact manifold, all vector fields on $N$ are complete: If $X$ is a vector field, then the flow curve $\gamma(t)$ defined by the initial value problem\[
 \gamma(0) = x \text{ and }
 \dot\gamma(t) = X_{\gamma(t)}
\]
exists for all $x\in N$ and $t\in \R$. Let us reprove a simple fact about Lyapunov functions on compact manifolds.
\begin{proposition}\label{lem:lyapunov}
 Let $f$ be a Lyapunov function on $N$ for $X \in C^\infty(N,TN)$, then the flow lines of $X$ start and end at singularities of $X$.
\end{proposition}

\begin{proof}
If $\gamma$ is a flow line of $X$ on $N$ with $\gamma(0)=x$ and $\gamma(t_0)$ a non-singular point,\[
\frac{d}{dt}(f\circ \gamma)|_{t=t_0} = Df_{\gamma(t_0)}(X) >0
\]
implies that $f\circ \gamma$ is increasing for all $t$. Since $N$ is compact and $f$ continuous, $\im (f \circ \gamma)$ is a bounded set. As $f\circ\gamma$ is always increasing, but bounded, we have\[
\lim_{t_0 \to \pm \infty} \frac{d}{dt}(f\circ \gamma)|_{t=t_0} = 0\;.
\]
By compactness we can furthermore find $x\in N$ and a sequence $t_n \in \R$ such that \[\lim_{n\to \infty} t_n = \infty\quad \text{and} \quad \lim_{n \to \infty} \gamma(t_n) = x \in N\;.\]
It follows that $x$ is singular, because $D f_x(X) = \lim_{n\to \infty}\frac{d}{dt}(f\circ \gamma)|_{t=t_n} =  0$. Since $X$ is smooth, we have $\|\dot \gamma(t)\| = \|X_{\gamma(t)}\|\to 0$ as $t\to\infty$, which implies \[\lim_{t\to\infty}\gamma(t) = x\;.\] A similar argument shows that $\gamma$ starts at a singularity.
\end{proof}

The following corollary from Proposition \ref{lem:lyapunov} describes the qualitative behavior of the flow of $D(\pi)(\nabla \vol)$ and $D(\pi)(\Psi (\nabla \vol))$. 
\begin{theorem}\label{ThmGradientFlowEnd}  The vector field $X =D(\pi) (\nabla \vol) = \nabla (\vol|_N)$ on $N$ is the gradient field of $f|_N$, and the flow lines of $X$ and $\tilde X = D(\pi)(\Psi (\nabla \vol))$ on $N \subset M$ start and end at singularities of $Y$.
\end{theorem}

\begin{proof}
We compute \[D(\pi) \nabla \vol_x = \nabla \vol_{x}-(\nabla \vol_{x}\cdot x)\cdot x \quad \text{for }x\in N.\] Since $x\in N \subset \R^{3k}$ is orthogonal to $X_x  \in T_x N \subset \R^{3k}$, we get\begin{align*}
\la \nabla (\vol|_N)_x, X_x \ra_N  & =  D\vol (X_x) = \la \nabla \vol_x,X_x \ra = \la D(\pi)(\nabla \vol_x) + (\nabla \vol_{x}\cdot x)\cdot x, X_x \ra\\
& = \la D(\pi)(\nabla \vol_x), X_x \ra  = \la D(\pi)(\nabla \vol_x), X_x \ra_\pi.
\end{align*}
This shows that $D(\pi)(\nabla \vol)$ is the gradient field $\nabla (\vol|_N)$ of $\vol|_N$. In particular, $\vol|_N$ is a Lyapunov function for $D(\pi)(\nabla \vol)$.

Therefore, we also have\begin{align*}
D (\vol|_N)_x \left(D(\pi)(\Psi (\nabla \vol_{\gamma(t)}))\right)& = \frac{1}{\sqrt{\nabla \vol_{\gamma(t)}}} D (\vol|_N)_x \left(D(\pi)(\nabla \vol_{\gamma(t)})\right)>0.
\end{align*}
Since $\Psi$ preserves singularities, $\vol|_N$ is also a Lyapunov function for $D(\pi)(\Psi(\nabla \vol))$. The theorem follows from Lemma \ref{lem:lyapunov}.
\end{proof}

We are interested in the normalization $\Psi(\nabla \vol)$ of $\nabla \vol$, because it is a scaling invariant vector field. One can also ask about the dynamic behavior of normalizations given by division of other powers of the norm and coordinate-wise normalizations. We expect them to have the same qualitative behavior. In the case of simply dividing by the norm of $X$, this dynamical behavior is described in Thom's gradient conjecture \cite{KMP}. The above theorem also shows that the flow of the vector fields $D(\pi)(\nabla\vol)$ and $D(\pi)(\Psi(\nabla\vol))$ optimize the quality measure $q=\vol\circ \pi$. See \cite{VartziotisHimpel2013b} how one can turn this observation into an efficient and global optimization-based smoothing method.

The vector fields $D(\pi)(\nabla\vol)$ and $D(\pi)(\Psi\nabla\vol)$ therefore have the same singularities and qualitative behavior as described Theorem~\ref{ThmGradientFlowEnd}. This suggests that the two vector fields are topologically equivalent, that is, there is a homeomorphism $N \to N$ carrying trajectories to trajectories and preserving the direction of increasing time. However, we do not see, how to prove such a strong result. Since the kernel of $D(\pi)_{x} \co \R^{3n} \to T_{x} N$ consists of the radial vectors, that is, vectors of the form $\lambda  x$, the singularities of these vector fields have the following characterization in terms of $\vol$.

\begin{lemma}\label{PropSing} The vector field
$D(\pi)(\vol)$ on $N \subset M$ has a singularity at $x\in N$ if and only if
\begin{equation}\label{EqSing}
\lambda x  = \nabla\vol_x \text{ for some constant } \lambda \in \R.
\end{equation}
\end{lemma}

\subsection{Symmetry}\label{SecSymmetry}

Mani \cite{M71} showed that for each polyhedral graph $G$, there exists a convex polyhedron $P_G$ such that every automorphism of $G$ is induced by an isometry of Euclidean space. We call $P_G$ a {\em symmetric} polyhedron. Examples for symmetric polyhedra include the platonic solids, semi-regular polyhedra, the right regular prism, right regular anti-prism and the regular pyramid. In this section every polyhedron $P$ is assumed to be symmetric.

Consider the {\em symmetry group} $\Gamma$ of a (symmetric) polyhedron $P \subset \R^3$ as a subgroup of isometry group of $\R^3$, that fixes $P$. Therefore, $g\in \Gamma$ acts on all points in $p \in {\R}^3$ as usual by $g\cdot p \in {\R}^3$. In particular, it naturally acts on a triangulation $T$ of $P$, and we can consider a $\Gamma$--invariant subset $\cT' \subset \cT$ of triangulations of polyhedra. In this way, a $\Gamma$--invariant triangulation of a polyhedron gives a simpler definition of the mean volume in Definition \ref{DefnMeanVol}\[
 \vol(x) = \frac{1}{|\cT'_x|}\sum_{T\in \cT'_x}\sum_{\tau \in T} \vol(\tau).
 \]
 Notice that the orbit of an arbitrary triangulation under $\Gamma$ yields a $\Gamma$--invariant subset of all triangulations.

Let us consider the two possible triangulations of the hexahedron consisting of 5 tetrahedra each and apply this theorem to this $\Gamma$--invariant subset $\cT$ of all triangulations.
\begin{equation}\label{EqGradHexa}
X^\cT_p = \displaystyle \frac{1}{12}\begin{pmatrix}
 \nu(2,5,4) + \nu(6,5,8,4,3,2)\\
 \nu(3,6,1) + \nu(7,6,5,1,4,3)\\
 \nu(4,7,2) + \nu(8,7,6,2,1,4) \\
 \nu(1,8,3) + \nu(5,8,7,3,2,1) \\
 \nu(1,6,8) + \nu(6,7,8,4,1,2)\\
 \nu(2,7,5) + \nu(7,8,5,1,2,3)\\
 \nu(3,8,6) + \nu(8,5,6,2,3,4)\\
 \nu(4,5,7) + \nu(5,6,7,3,4,1)
\end{pmatrix} \quad \text{for 
\raisebox{-1cm}{\def\svgwidth{3cm}
%% Creator: Inkscape inkscape 0.48.3.1, www.inkscape.org
%% Accompanies image file 'hexa.pdf' (pdf, eps, ps)
\begingroup%
  \makeatletter%
  \providecommand\color[2][]{%
    \errmessage{(Inkscape) Color is used for the text in Inkscape, but the package 'color.sty' is not loaded}%
    \renewcommand\color[2][]{}%
  }%
  \providecommand\transparent[1]{%
    \errmessage{(Inkscape) Transparency is used (non-zero) for the text in Inkscape, but the package 'transparent.sty' is not loaded}%
    \renewcommand\transparent[1]{}%
  }%
  \providecommand\rotatebox[2]{#2}%
  \ifx\svgwidth\undefined%
    \setlength{\unitlength}{123.92514648bp}%
    \ifx\svgscale\undefined%
      \relax%
    \else%
      \setlength{\unitlength}{\unitlength * \real{\svgscale}}%
    \fi%
  \else%
    \setlength{\unitlength}{\svgwidth}%
  \fi%
  \global\let\svgwidth\undefined%
  \global\let\svgscale\undefined%
  \makeatother%
  \begin{picture}(1,0.78723177)%
    \put(0,0){\includegraphics[width=\unitlength]{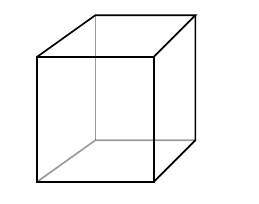}}%
    \put(0.079027,0.0182696){\color[rgb]{0,0,0}\makebox(0,0)[lb]{\smash{$x_1$}}}%
    \put(0.54698034,0.0182696){\color[rgb]{0,0,0}\makebox(0,0)[lb]{\smash{$x_2$}}}%
    \put(0.75685553,0.17965734){\color[rgb]{0,0,0}\makebox(0,0)[lb]{\smash{$x_3$}}}%
    \put(0.36952494,0.17965734){\color[rgb]{0,0,0}\makebox(0,0)[lb]{\smash{$\cg x_4$}}}%
    \put(-0.00643029,0.53471038){\color[rgb]{0,0,0}\makebox(0,0)[lb]{\smash{$x_5$}}}%
    \put(0.60399935,0.50812041){\color[rgb]{0,0,0}\makebox(0,0)[lb]{\smash{$x_6$}}}%
    \put(0.77391846,0.70178571){\color[rgb]{0,0,0}\makebox(0,0)[lb]{\smash{$x_7$}}}%
    \put(0.2081372,0.72837568){\color[rgb]{0,0,0}\makebox(0,0)[lb]{\smash{$x_8$}}}%
  \end{picture}%
\endgroup%
}.}
\end{equation}

In addition to the hexahedron, it is now a straight-forward task to write down and check gradient fields for the rest of the platonic solids. In \S\ref{SecIcoDode} we mention a gradient field for the icosahedron, which is a variation of the gradient field of the mean volume. We leave the explicit gradient field for the dodecahedron to the interested reader, because the expression is considerably longer.

\subsection{Smoothing of hybrid meshes}

We have constructed elementary transformations for tetrahedra, hexahedra, pyramids, prisms and octahedra. We will see in \S\ref{SecRegularity} and \S\ref{SecDyna}, that they make individual volume elements more regular. Moreover, we will single out the vector fields, which exhibit the fastest convergence in numerical tests for each of the four main volume elements. Using the corresponding GETMe transformations, we get a smoothing algorithm for hybrid meshes as described in \S\ref{SecGETme}, whose performance we compare to the original results from GETMe transformations in \S\ref{SecDyna}.

\section{Dynamics}\label{SecDyna}

So far, we have constructed a quality measure based on the mean volume of each volume element, related its gradient to the tetrahedral GETMe smoothing algorithm and generalized it to hybrid meshes. We are left with analyzing the dynamics, both theoretically and numerically. For a complete theoretical treatment, we will need to compute the singularities of the vector field. In view of the discussion on symmetry, it is not surprising that regular polyhedra are singularities of the vector field. In fact, numerical tests in \S\ref{Sec:NumericalTests} show that regular polyhedra are sinks of its flow. However, not every singularity is a regular polyhedron. Even though this discussion is elementary, the geometric arguments are rather intricate. Due to the technical nature, we postpone this topic to \S\ref{SecRegularity}. We start by showing that the signed volume for tetrahedral volume elements is a Morse-Bott function on the 8--sphere $N$, that is, a smooth function on $N$ whose critical set is a closed 
submanifold and whose Hessian is non-degenerate in the normal direction. In particular, this explains the untangling and regularizing behavior of its flow. We will then see, what we can say about the dynamics for other volume elements from a theoretical point of view, before we provide numerical evidence.

\subsection{Volume as a Morse-Bott function for tetrahedra}\label{SecMorse}

In general, the group $\SO(3)$ acts on $(x_1,\ldots,x_n) \in M$ by rotating all $x_i$ simultaneously, and the vector field $\nabla\vol$ is invariant under rotation. Due to the rotational invariance, all critical points of the mean volume on $N$ are therefore degenerate. The action of $\SO(3)$ on the  $(3n-4)$--sphere $N$ is proper and free as long as the $x_i$ are not collinear: The stabilizer $S_x = \{ A \in \SO(3) \mid x\cdot A = x\}$ for every collinear $x\in M$ is $S^1$. The collinear tetrahedra consist of the $n$--dimensional vector space $V$ of $n$ ordered points on each line through the origin in $\R^3$. These rays can be parameterized by $S^2$. Since tetrahedra with all vertices having the same coordinates are excluded in $M$, this parametrization by $S^2$ is one-to-one. The map $\pi\co M \to N$ first translates the $n$ points on the ray, so that their centroid is the origin, then normalizes the result to have length $1$. This corresponds to a surjection of $V$ onto an $n-1$--dimensional sub-vector 
space $W$ 
of $V$, followed by a projection to the $n-2$--dimensional unit-sphere in $W$. Therefore, the subspace  $S_c$ of collinear $x \in N$ is diffeomorphic to $S^2 \times S^{n-2}$ via $\pi$. Clearly, the set of singular tetrahedra with vanishing volume consists of the collinear ones and is diffeomorphic to $S^2 \times S^2$.  In general, $S_c$ is only a subset 
of the critical points of the mean volume. It is not difficult to see, that the negatively and the positively oriented regular tetrahedron is also a singularity. Since the tetrahedra can be rotated, these two families are parametrized by $S^2$.

In order to show that the (mean) volume on $N$ is Morse-Bott for tetrahedra, we need to compute the eigenvalues of the Hessian. The Hessian of the volume on $M$ is equal to the gradient of $\nabla \vol$.  In order to compute the eigenvalues of the Hessian $\Hess(\vol|_N)$ of the restriction of $\vol\co M \to \R$ to $N$ at a regular tetrahedral volume element $x$, we can compute the eigenvalues of the Hessian of $\vol\circ \pi$ and keep in mind, that six eigenvalues are zero; three from the translation invariance, two from the rotation invariance and one from the scaling invariance. On the level of vector fields we need to take the differential $D(\pi)(X_x)$ of $X_x$ via $D(\pi) \co T_x M \to T_{\pi(x)} N$ before we take the gradient. 
In order to compute the eigenvalues of $\Hess(\vol|_N)$ at positively oriented regular polyhedra, we can therefore compute the eigenvalues of the Jacobian matrix given by the gradient of $D(\pi)(X_x) \subset T_x M$. The result of this computation is summarized in Table~\ref{TabEigen}. By symmetry, the eigenvalues for the negatively oriented polyhedra are the same, only with the opposite sign. This agrees with the observation that the mean volume has a global minimum at the negatively oriented regular polyhedron and a global maximum at the positively oriented regular polyhedron.

\begin{table}[ht]
\centering
\setlength\extrarowheight{0.5ex}
\begin{tabular}{l|c|c|l}
Type & X & $n$ & $(\mu,m) = $ eigenvalue $\mu$ of the $\Hess(\vol|_N)$ with multiplicity $m$\\
\hline
tetrahedron & \eqref{EqTetraX} & $4 \cdot 3$ & $(-\frac{2}{\sqrt{27}},6)$\\
pyramid & \eqref{EqGradPyra} & $5 \cdot 3$ & $(- \sqrt{\frac{5}{27}},6), (- \frac{\sqrt{5}}{6\cdot\sqrt{3}},3)$\\
octahedron &  \eqref{EqGradOcta} & $6 \cdot 3$ & $(-\sqrt{\frac{8}{27}},6), (-\sqrt{\frac{2}{27}},6)$\\
prism & \eqref{EqGradPrism} &  $6 \cdot 3$ & $(-\sqrt{\frac{2}{27}},6), (-\frac{\sqrt{8}}{9\cdot\sqrt{3}},2),(-\frac{1}{3\cdot\sqrt{6}},2),(-\frac{\sqrt{2}}{9\cdot\sqrt{3}},2)$\\
hexahedron & \eqref{EqGradHexa} & $8 \cdot 3$ & $(-\frac{1}{\sqrt{6}},6), (-\frac{5}{6\cdot\sqrt{6}},1), (-\sqrt{\frac{2}{27}},3),(-\frac{1}{\sqrt{24}},3),(-\frac{1}{\sqrt{56}},5)$\\
hexahedron & \eqref{EqGradHexa2} & $8 \cdot 3$ & $(-\sqrt{\frac{32}{3}},6), (-\sqrt{\frac{8}{3}},12)$
\end{tabular}\\[2ex]
\caption{The eigenvalues of the Hessian at positively oriented regular polyhedra}
\label{TabEigen}
\end{table}

We also compute that the Hessian at collinear tetrahedra has two positive and two negative eigenvalues.
This shows that the (mean) volume is a Morse-Bott function on the $8$--sphere $N$, more precisely, the critical sets are submanifolds $S^2$, $S^2\times S^2$ and $S^2$ at the maximum, level 0 and the minimum, respectively, and the Hessian is non-degenerate in the normal direction. Furthermore, the positively oriented regular tetrahedron is a sink of the gradient flow as well as of the flow of $\Psi(\nabla \vol)$, while the stable and unstable manifold are two-dimensional for collinear tetrahedra. Generically, almost all tetrahedra are therefore regularized by the volume gradient flow, more precisely, the set of all tetrahedra, which are regularized, is open and dense in $N$. Alternatively, a small perturbation of the mean volume in a neighborhood of the level 0 singularities is a Morse function with two singularities consisting of the positively and negatively oriented tetrahedra. Unfortunately, for other types of polyhedra the mean volume function is not Morse-Bott in general, despite it being constructed 
from 
the 
tetrahedral volume.

\subsection{Other volume elements}\label{Sec:DynaOther}

It is easy to check that regular polyhedra are singularities of the gradient of the mean volume in the sense of Theorem~\ref{ThmGradientFlowEnd}. Furthermore, Table~\ref{TabEigen} shows, that regular tetrahedra, pyramids, octahedra, prisms and hexahedra are sinks (and sources) of the corresponding vector fields. However, not every singularity is a regular polyhedron. The rigorous identification of the singularities is a difficult problem and the solution is bound to be technical. For this purpose, the notion of regularity needs to be discussed in more detail and requires a definition. Therefore, we postpone this topic to \S\ref{SecRegularity} and first present numerical tests in \S\ref{Sec:NumericalTests} confirming the regularizing behavior of the gradient of the mean volume. We show rigorously in \S\ref{SecRegularity}, that there is only one singularity with positive mean volume in the case of the tetrahedron, the pyramid, the octahedron, the prism and the hexahedron, namely the regular one. Then, this 
singularity is a sink, while the regular representative with negative mean volume is a source. It follows, that every volume element with positive mean volume will be regularized.

In order to have a complete picture about the dynamics and a theoretical discussion about the untangling effect of the gradient flow, we need to address the singularities with vanishing mean volume, that is, the level 0 singularities, and their stable and unstable sets. Clearly, there are a lot of (necessarily invalid) volume elements, which would flow towards level 0 singularities. However, we expect to be able to regularize most arbitrarily small perturbations of the initial volume element. This is a heuristic statement. In order to make this mathematically precise, we would have to discuss measure theory or probability theory, which would allow us to say, that the set of volume elements, which is regularized has measure one, or that the volume elements will be regularized with probability one. However, the first necessity is the computation of the level 0 singularities. We expect these singularities to be unions of submanifolds of $N$, and that the stable sets are also union of submanifolds. We expect the 
maximal dimension of the stable manifolds to be 
less 
than 
the dimension of $N$. In particular, the complement of these stable manifolds would then be open and dense in $N$, which suggests, that it has measure one. This is numerically confirmed by our tests in \S\ref{Sec:NumericalTests}. Unfortunately, this is not so easy to prove rigorously. Already the computation of the level 0 singularities is difficult.

\begin{figure}[ht]
\def\svgwidth{7cm}
\begin{center}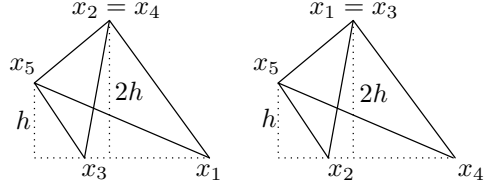\end{center}\caption{A family of (planar) singularities for pyramids.\label{pyra_sing}}
\end{figure}

Let us demonstrate the computation of the level 0 singularities for the pyramid. Consider $\nabla\vol_x = 0$ in the case of the pyramid, where $x$ is not collinear. If $\nu(1,2,3,4) = 0$ then $x_2=x_4$ or $x_1=x_3$. Suppose that $x_2=x_4$ and $\nu(5,1,3) \neq 0$,  then by $0 = \nu(5,1,3) + \nu (5,1,4,3) = 2\nu(5,1,3) + \nu(1,4,3)$, so that $x_2=x_4$ lies in the plane spanned by $x_1$, $x_3$ and $x_5$, such that the height of the triangle $C(1,4,3)$ is twice the height of the triangle $C(1,5,3)$. Therefore, the two planar constellations given in Figure~\ref{pyra_sing} are singular pyramids. The case $\nu(1,2,3,4) \neq 0$ is not possible for level 0 singularities, which is proven in more generality in Lemma \ref{Lem:OptEqiv}. The submanifold of $N$ where $\nu(1,2,3,4)=0$ and $x_1=x_3$ is diffeomorphic to the Cartesian product of real projective 3--space with the 3--sphere $\SO(3) \times S^3 = \R P^3 \times S^3$, but it intersects other submanifolds of the singular set. Therefore, the singular set is a union of 
three submanifolds of $N$, namely $S^2 \times S^3$ and two 
copies of $\SO(3) \times S^3$.

In applications, however, we usually do not start with volume elements of negative mean volume. In other words, if an initial mesh has volume elements negative mean volume, then we expect, that there is a problem with the setup. If we want to apply the flow to individual valid volume elements, they will have positive mean volume and can be regularized. In view of the theoretical proof in \S\ref{SecRegularity}, that the only singular tetrahedra, pyramids and octahedra with positive volume are regular polyhedra, the above discussion can be summarized in the following mathematically rigorous result.

\begin{theorem} The mean volume gradient flow regularizes volume elements based on tetrahedra, pyramids and octahedra, as long as their volume is positive. The gradient flow of the gradient vector $Y$ given in \eqref{EqGradHexa2} regularizes all hexahedral volume elements $x\in N$ if they satisfy $f^Y >0$, where $f^Y$ is defined in \eqref{EqHexa2}.
\end{theorem}

\subsection{Numerical tests}\label{Sec:NumericalTests}

In practice it is not important, that the vector field is actually a gradient field, therefore we have changed the vector field for the prism, the hexahedron and the pyramid slightly so that we observe faster convergence towards a regular representative. In the case of a regular pyramid and prism, we also decided that the ratio of the edges should be one. We have furthermore reduced the number of necessary cross products by making use of the centroid preserving behavior and Remark~\ref{RemIdentities}. Only the vector field for the prism is not centroid preserving. We therefore base our tests on the vector fields defined in \eqref{eq:X_tet}, \eqref{eq:X_pyr}, \eqref{eq:X_pri} and \eqref{eq:X_hex}.   Each vector field $X=(X_1,\dots,X_k)$ is first multiplied by $1/\sqrt{\|X\|}$ to make it scaling-invariant, then its growth is harmonized by multiplying $X/\sqrt{\|X\|}$ by a scalar $C$ depending on the vector field, so that $|x'_i-x'_j|  = 2 |x_i-x_j|$ for a regular polyhedron $x$ and $\sigma = 1$, where
 \[
   x' =  x + \sigma\frac{C}{\sqrt{\|X\|}} X.
 \]

In the case of the tetrahedron we only need 3 cross products in the vector field
\begin{equation}\label{eq:X_tet}
X^\text{tet}_x = \begin{pmatrix}
\nu(4,3,2)\\
\nu(4,1,3)\\
\nu(4,2,1)\\
\nu(1,2,3)
\end{pmatrix} = \begin{pmatrix}
\nu(4,3,2)\\
\nu(4,1,3)\\
\nu(4,2,1)\\
-\nu(4,3,2)-\nu(4,1,3)-\nu(4,2,1)\\
\end{pmatrix}
\end{equation}

In the case of the pyramid only 5 cross products are necessary
\begin{equation}\label{eq:X_pyr}
X^\text{pyr}_x =
\begin{pmatrix}
\nu(5,4,3,2)+4\cdot \nu(5,4,2)\\
\nu(5,1,4,3)+4\cdot \nu(5,1,3)\\
\nu(5,2,1,4)-4\cdot \nu(5,4,2)\\
\nu(5,3,2,1)-4\cdot \nu(5,1,3) \\
2 \cdot\nu(1,2,3,4)\end{pmatrix}
=
\begin{pmatrix} 
\nu(5,4,3,2)+4\cdot \nu(5,4,2)\\
\nu(5,1,4,3)+4\cdot \nu(5,1,3)\\
-\nu(1,2,3,4)-\nu(5,4,3,2)-4\cdot \nu(5,4,2)\\
-\nu(1,2,3,4)-\nu(5,1,4,3)-4\cdot \nu(5,1,3) \\
2 \cdot\nu(1,2,3,4)\end{pmatrix}
\end{equation}

For the prism 5 cross products suffice
\begin{equation}\label{eq:X_pri}
X^\text{pri}_x = \begin{pmatrix}
\nu(2,5,6,3) + 3\cdot\nu(5,4,6)\\
 \nu(3,6,4,1)+ 3\cdot\nu(5,4,6)\\
 \nu(1,4,5,2)+ 3\cdot\nu(5,4,6)\\
 \nu(2,5,6,3)+ 3\cdot\nu(1,2,3)\\
 \nu(3,6,4,1)+ 3\cdot\nu(1,2,3)\\
 \nu(1,4,5,2)+ 3\cdot\nu(1,2,3)\\
\end{pmatrix}\\
\end{equation}

In the case of the hexahedron we only need to compute 5 cross products. Given \begin{align*}
X_1 =  \nu(1,5,8,4)+\nu(1,2,6,5)+\nu(1,4,3,2),\\
X_2 = \nu(1,2,6,5)+\nu(2,3,7,6)+\nu(1,4,3,2),\\
X_3 = \nu(2,3,7,6)+\nu(3,4,8,7)+\nu(1,4,3,2),\\
X_4 = \nu(3,4,8,7)+\nu(1,5,8,4)+\nu(1,4,3,2),
\end{align*}
we have
\begin{equation}\label{eq:X_hex}
X^\text{hex}_x = (X_1,X_2,X_3,X_4,-X_3,-X_4,-X_1,-X_3)^t.
\end{equation}

For the above transformations the convergence has
been tested analogously to the approach given in \cite{VartziotisWipper2011Mixed}: The tests
are based on generating 100,000 random initial elements $E_j$ of each type and
taking 101 equidistant values for the scaling factor $\sigma\in[0,1]$.  For
each pair $(E_j,\sigma_k)$ taken from the Cartesian product of all elements
and all $\sigma$-values, the geometric transformation using the scaling factor
$\sigma_k$ has been iteratively applied starting with the initial element
$E_j$ until the mean ratio quality number of the resulting element deviated
less than $10^{-6}$ from the ideal value one or the number of iterations
exceeded 100. For each transformation and the associated element type this
test resulted in 10,100,000 transformation cycles, each represented by its
number of iterations. The validity of the initial element is not a necessary condition for regularization.
Arbitrary random nodes lead to a large amount of
invalid elements (99.6\% in the case of hexahedral elements). However, we get an excellent convergence of all of these elements comparable to \cite[Figure~2]{VartziotisWipperSchwald2009} and \cite[Figure~4]{VartziotisWipper2011Mixed}. In particular, this shows that the transformation untangles individual elements. 

\begin{figure}
\def\svgwidth{\textwidth}
\centering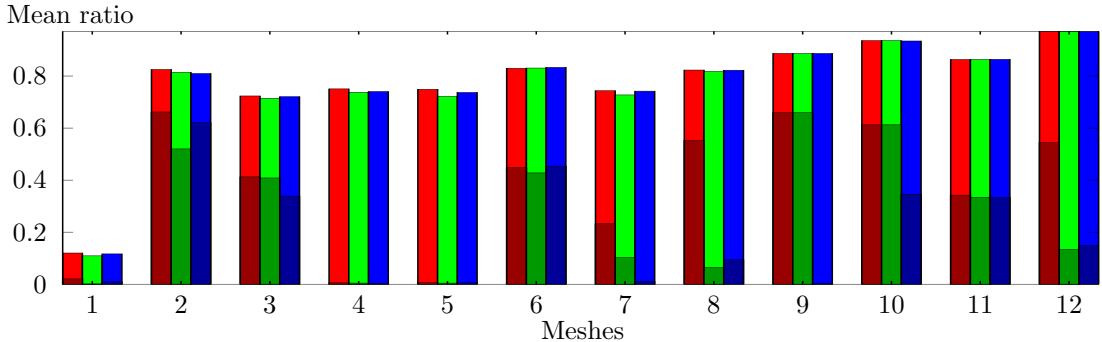
\captionsetup{singlelinecheck=off}\caption[Caption mit itemize]{Comparison of the arithmetic mean and the minimum of all mean ratio numbers for 12 test cases using the smoothing in the original paper (left), the dual element-based smoothing with the parameters from \cite{VartziotisWipperPapadrakakis2013} (middle) and the new GETMe smoothing based on the mean volume (right), where  \begin{itemize}
   \item Meshes 1--5 correspond to Figures 5, 8, 12, 18, 19 in \cite{VartziotisWipper2011Hex},
   \item Meshes 6--8 correspond to Figures 6, 10, 14 in \cite{VartziotisWipper2011Mixed},
   \item Meshes 9--10 correspond to the hexahedral and tetrahedral meshes in \cite[Figure~10]{VartziotisWipperPapadrakakis2013}, and
   \item Meshes 11-12 correspond to Figures 5, 8 in \cite{VartziotisPapadrakakis2013}.\end{itemize}\label{fig:GETMe_Comparison}}\end{figure}

In order to test the performance of the resulting simultaneous and sequential GETMe smoothing for hybrid meshes, we used all 32 test cases and results from \cite{VartziotisWipperSchwald2009,VartziotisWipper2011Hex,VartziotisWipper2011Mixed,VartziotisWipperPapadrakakis2013,VartziotisPapadrakakis2013} and compared the performance of our new smoothing algorithm for $\sigma =1$ with the dual element-based GETMe algorithm using the same (optimized) parameters as in \cite{VartziotisWipperPapadrakakis2013}. Some of the results are shown in Figure~\ref{fig:GETMe_Comparison}. The run-times are comparable, therefore only the arithmetic mean and the minimum of all mean ratio numbers are given. Just like in all previous papers on GETMe smoothings, the boundary nodes have been 
kept fixed, which explains why some examples only show poor improvements. We see, that the dual element-based and our new smoothing 
algorithm produce meshes of a comparable mean ratio quality in about the same amount of time. Only Mesh \#9 has elements with exceptionally low quality, but for this example the parameters had been manually optimized for the dual-element based transformation as in \cite{VartziotisWipperPapadrakakis2013}. Still, the arithmetic means are comparable throughout the examples. We should however note, that our algorithm is based on a different quality measure, therefore the mean ratio quality measure applied to our results is sometimes not as good as one could hope for.

\section{Theoretical analysis of regularity} \label{SecRegularity}

In the theoretical treatment of the dynamics for the vector fields given in \eqref{eq:X_tet}, \eqref{eq:X_pyr}, \eqref{eq:X_pri} and \eqref{eq:X_hex} we owe a rigorous computation of the singularities. We have seen in \S\ref{Sec:DynaOther}, that these sets become very complicated. Therefore, we will focus on singularities for polyhedra with positive volume, which include valid polyhedra. We will show, that singularities of positive volume are regular polyhedra.

\subsection{Terminology}

Intuitively, symmetry and regularity are related notions. One definition found in the literature is, that a polyhedron is regular, if the action of its symmetry group acts transitively on its flags, so that the platonic solids are the only regular, convex polyhedra \cite{E}. A flag is a connected set of elements of each dimension---in the case of a polyhedron it consists of the body, a face, an edge of the face, a vertex of the edge, and the null polytope. Nevertheless, prisms and pyramids are also sometimes called regular, if they are symmetric. This notion can also be extended to isogonal polytopes, that is, polytopes, for which their symmetry group acts transitively on its vertices. Let us therefore call a pyramid, a prism and an isogonal polyhedron regular, if it is symmetric (and convex). It has been observed heuristically in \cite{V}, that the iterative application of certain actions induced by the symmetry group on polygons regularizes them. In a way this is our hope for the flow of the gradient 
fields of the 
mean volume for a $\Gamma$--invariant set of triangulations $\cT$.

In the case of tetrahedra, we will show directly in \S\ref{SecTetra} that the non-collinear singularities are regular tetrahedra. In particular, we see by Table~\ref{TabEigen} that the positively oriented one is a sink, and the negatively oriented one is a source of the gradient flow. It is tempting to think that this generalizes to other polyhedra with non-vanishing mean volume, namely that the singularities of the mean volume flow on $N$ are regular polyhedra in some sense. Unfortunately, this naive approach will only work for rather simple polyhedra, even if they have a large symmetry group. For example, in the case of the icosahedron, there are a lot more singularities of the mean volume function than just the regular icosahedron, and, in particular, the regular icosahedron does not maximize the mean volume function on $N$. In this case, a slightly more complicated gradient field seems to have the desired symmetric icosahedra as singularities: the regular icosahedron and the great icosahedron (one of the 
Kepler–Poinsot polyhedra).

While we have confirmed numerically, that the gradient flow of the mean volume exhibits a regularizing dynamic behavior for the tetrahedron, the pyramid, the octahedron, the hexahedron and even the dodecahedron, we can rigorously confirm this behavior only for the tetrahedron, the pyramid and the octahedron. For the hexahedron, we can show it for a slightly different gradient field. Surprisingly and regrettably, we have not yet been able to prove it for the prism, even though numerical tests, in particular a computation of the eigenvalues of the Hessian of the collinear prisms, suggest, that the mean volume is a Morse-Bott function just like we prove it for tetrahedra in \S\ref{SecMorse}. The icosahedron and the dodecahedron have only been studied numerically.

Based on Lemma~\ref{PropSing} we make the following definition:

\begin{definition}\label{DefnReg} Let $X$ be a vector field on $N = M/\sim$ given by a linear combination of cross products. The volume element $x$ is {\em $X$--optimal}, if there exists $y \sim x$ satisfying
\begin{equation}\label{EqCond}
\lambda y = X_y \quad \text{for some }\lambda \neq 0.
\end{equation}
In particular, no collinear $x$ is $X$--optimal. The sign of $\lambda$ determines the orientation of $x$.
\end{definition}

The following lemma shows, that the $X$--optimal volume elements $x$ are the most interesting singularities of the gradient of the mean volume.

\begin{lemma}\label{Lem:OptEqiv}
Consider a singularity $x$ of $D(\pi)(\Psi\nabla\vol)$. Then $x$ has non-zero mean volume, if and only if it is $\nabla\vol$--optimal.
%If a volume element $x$ is a singularity of $D(\pi)(\Psi\nabla\vol)$ with non-zero mean volume, then it is $\nabla\vol$--optimal.
\end{lemma}

\begin{proof}
 If $x$ has non-zero mean volume, then rescaling changes the volume. In particular, the tangent vector $X\in T_xM$ corresponding to shifting the nodes satisfies $\la\nabla\vol_x,X_x\ra = D(\vol_x)(X_x) \neq 0$. Therefore, $\nabla\vol_x \neq 0$. Lemma~\ref{PropSing} and Theorem~\ref{ThmGradientFlowEnd} then show that $x$ is $\nabla\vol$--optimal. On the other hand, if $x$ has zero mean volume, rescaling will preserve the volume. Therefore $\nabla\vol_x = 0$. This finishes the proof.
\end{proof}

We expect that the $X$--optimal volume elements are sinks or sources depending on the sign of $\lambda$. It is difficult in general to determine the shape of the singularities. Let us therefore discuss $X$--optimality for the platonic solids, the prism and the pyramid, when $X$ is the gradient of the mean volume or a slight variation of it. In order to show, that an $X$--optimal volume element is regular, we first deduce its symmetry from elementary geometric arguments, which can be technical and tricky. We will only provide some of the proofs to show the nature of the geometric arguments and leave a few of the details in other proofs to the interested reader.

\subsection{The tetrahedron}\label{SecTetra}

The tetrahedron is the simplest case. We show in \S\ref{SecMorse}, that the mean volume is a Morse-Bott function on the sphere $N$. In particular, the dynamics of its gradient flow are well-understood, see \S\ref{SecDyna}.

\begin{theorem} The non-collinear singularities of the gradient of the (mean) volume for tetrahedra are regular tetrahedra.
\end{theorem}
\begin{proof}
Consider the gradient  $\nabla \vol$ of the volume given in \eqref{EqTetraX}. If $x$ is a non-collinear singularity, it is $(\nabla \vol)$--optimal. Let us assume without loss of generality that $x_4=0$. Then\[
\frac{6\lambda}{2} (x_1-x_2) = \frac{1}{2} (\nu(4,3,2) - \nu(4,1,3))\\
 = x_3 \times \left(\frac{1}{2}(x_1+x_2)\right).
\]
Therefore, $x_1-x_2$ is orthogonal to the hyperplane spanned by $x_3$ and $\frac{1}{2}(x_1+x_2)$. In particular, the face with vertices $x_1$, $x_2$ and $x_3$ is isosceles, see Figure~\ref{tetrarect}. By symmetry, all edges of the tetrahedron are of equal length, and therefore the tetrahedron is regular.\begin{figure}[ht]
\centering\def\svgwidth{2.2cm}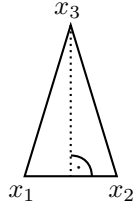\caption{Isosceles triangle}\label{tetrarect}
\end{figure}
\end{proof}

\subsection{The pyramid}\label{SecPyra}

Consider the gradient $(\nabla \vol)$ of the mean volume function for pyramids given in \eqref{EqGradPyra}. Suppose $x$ is $\nabla \vol$--optimal. We compute $x_4-x_3 = x_2-x_1$. Therefore, the base of a singularity is a (planar) parallelogram. If we translate $x$ so that $x_5=0$, a geometric argument as before shows that\begin{align*}
12\lambda(x_4-x_2) &=2 \nu(5,3,1) +\nu(5,3,2,1) + \nu(5,3,4,1)\\
&= 2x_3 \times x_1 + (x_3-x_1) \times (x_2+x_4)\\
& = (x_3-x_1) \times (x_1 + x_3 + x_2 + x_4) .
\end{align*}
The right-hand side is orthogonal to $x_3-x_1$. Therefore, the base is a square, and $x$ is a symmetric pyramid, whose height is determined by the following result.

\begin{theorem}\label{ThmRegPyra} Every $(\nabla\vol)$--optimal pyramid is of the form $\pi(x) \in N$, where\[
x = \begin{pmatrix} (0,0,0)\\
(2,0,0)\\
(2,2,0)\\
(0,2,0)\\
(1,1,\sqrt{5})
\end{pmatrix}.
\]
\end{theorem}

{\em Proof}.  Consider a pyramid of height $h$ with square base of area $a^2$ and apex orthogonally above the middle of the base. If $x$ is $(\nabla\vol)$--optimal, then we compute $\|\nu(1,2,3,4) \|= 2a^2$ and $\|\nu(5,2,4)\|= h\cdot \|x_4-x_2\|$. For $\lambda \neq 0$ as in Definition \ref{DefnReg} we have\begin{align*}
24\lambda h & = 24\lambda\left\|x_5 - \frac{x_1+x_2+x_3+x_4}{4}\right\| = \|2\nu(1,2,3,4) - \frac{1}{2}\nu(4,3,2,1)\| = 5 a^2\\
\tag*{and} 24\lambda \|x_3-x_1 \| & = \|4\nu(5,2,4) + \nu(2,1,4)+\nu(3,4,2) \| = 4h \|x_4-x_2\|.
\end{align*}
Since $\|x_3-x_1\| = \|x_4-x_2\|$, we compute \[\frac{h}{a} =\frac{\sqrt{5}}{2}.\qquad\endproof%\qedhere
\]

\begin{rem}\label{RemPyra}{\em 
From the perspective of the finite element method, it is desirable, that all edges have the same length. We have two options to take care of this issue. We can ignore it entirely by arguing, that the edges connected with the apex for $x$ as in Theorem \ref{ThmRegPyra} have length $\sqrt{7} \approx 2.64575$, which might be close enough to 2 for mesh smoothing applications. Alternatively, if we are willing to let go of the theoretical advantages of a gradient field, it is straight-forward to vary the vectors in the gradient field so that all edges have the same length and the optimal singularities are preserved. Clearly, we have a lot of (topologically) equivalent transformations at our disposal and we can control the shape of the optimal pyramid, if we could only prove topological equivalence of the corresponding vector fields.}
\end{rem}

\subsection{The octahedron}\label{SecOcta}

Consider the gradient $\nabla\vol$ of the mean volume for the octahedron given in Equation~\eqref{EqGradOcta}. Let $x$ be a $(\nabla\vol)$--optimal volume element. Then, we can readily calculate that $x_3-x_6 = x_1-x_5$. Together with analogous statements for the other edges, we see that $C(2,3,4,5)$ and the other equator curves are a parallelograms. Furthermore, $x_1-x_6 = 2 \nu(2,3,4,5)$, therefore $x_1-x_6$ is orthogonal to the plane, in which $C(2,3,4,5)$ lies. Together with the analogous statements for the other cases this implies orthogonality of the planes, in which $C(2,3,4,5)$, $C(1,2,6,4)$ and $C(1,3,6,5)$ lie.

To prove the regularity of the octahedron, it remains to check that the area of the rhombi circumscribed by $C(2,3,4,5)$, $C(1,2,6,4)$ and $C(1,3,6,5)$ are all equal. We compute\begin{align*}
\|\nu(1,3,6,5)\| & = \frac{1}{2}\|x_1-x_6\| \cdot \|x_3-x_5\| = \frac{(6\lambda)^{-2}}{2}\| 2 \nu(2,3,4,5)\| \cdot \| \nu (1,2,6,4) - \nu(1,4,6,2) \|\\
&  = 2(6\lambda)^{-2} \|\nu(2,3,4,5)\| \cdot \|\nu(1,2,6,4)\|.
\end{align*}
By symmetry we also have $\|\nu(2,3,4,5)\| = 2(6\lambda)^{-2} \|\nu(1,3,6,5)\|\cdot \|\nu(1,2,6,4)\|$. Therefore, once again by symmetry,\[
\|\nu(2,3,4,5)\| = \|\nu(1,3,6,5)\| =  \|\nu(1,2,6,4)\| = \frac{1}{2}(6\lambda)^{2}.
\]
This shows the following.

\begin{theorem}
Every $(\nabla\vol)$--optimal $x$ is a regular octahedron. 
\end{theorem}

\subsection{The prism}\label{SecPrism}

The prism turns out to be the most difficult case. We have not been able to rigorously compute the singularities of the gradient $\nabla \vol$ of the mean volume on $N$ given in Equation \eqref{EqGradPrism}, even though the numerical tests confirm that the $\nabla\vol$--optimal prisms are regular. However, the optimal singularities for the simpler vector field
\begin{equation}\label{EqGradPrism2}
Y_x = \begin{pmatrix}
\nu(3,2,5,4,6)\\
\nu(1,3,6,5,4) \\
\nu(2,1,4,6,5)\\
\nu(5,6,3,1,2)\\
\nu(6,4,1,2,3)\\
\nu(4,5,2,3,1)\\
\end{pmatrix}
\end{equation}
can be easily computed to be regular prisms. We leave the proof to the interested reader.

Due to the numerical tests in \S\ref{Sec:NumericalTests} we expect the $\nabla\vol$--optimal prism to be symmetric. Just like in the case of the pyramid, we can compute that for a symmetric $\nabla\vol$--optimal prism we have $h/a = \sqrt{2/3}$, where $h$ is the height of the prism and $a$ is the length of each side of the base.

\begin{rem}\label{RemPrism}{\em 
As in the case of the pyramid, a prism with all edges of the same length might be more desirable. It is straight-forward to construct a non-gradient field from the gradient field, such that the optimal singularities have the desired edge-length.}
\end{rem}

\subsection{The hexahedron}\label{hexa}

We have not been able to compute the optimal polyhedra for the gradient of the mean volume given in Equation \eqref{EqGradHexa}. However, the vector field  
\begin{equation}\label{EqGradHexa2}
Y_x = \frac{1}{2}
\begin{pmatrix}
\nu(3,6,8) + \nu(2,5,4) + \nu(6,5,8,4,3,2)\\
\nu(4,7,5) + \nu(3,6,1) + \nu(7,6,5,1,4,3)\\
\nu(1,8,6)+ \nu(4,7,2) + \nu(8,7,6,2,1,4) \\
\nu(2,5,7)+ \nu(1,8,3) + \nu(5,8,7,3,2,1) \\
\nu(2,7,4) + \nu(1,6,8) + \nu(6,7,8,4,1,2)\\
\nu(3,8,1)+ \nu(2,7,5) + \nu(7,8,5,1,2,3)\\
\nu(4,5,2) + \nu(3,8,6) + \nu(8,5,6,2,3,4)\\
\nu(1,6,3)+ \nu(4,5,7) + \nu(5,6,7,3,4,1)
\end{pmatrix}
\end{equation}
differs from the gradient of the mean volume only slightly and is still a gradient field: If $\tau_1(x) = (x_3,x_6,x_8,x_1)$ and $\tau_2(x) = (x_2,x_4,x_5,x_7)$, then we have
\begin{equation}\label{EqHexa2}
Y_x = 6 \nabla f^Y_x \quad\text{for }f^Y \equiv \vol + \frac{1}{2} \vol\circ\tau_1 + \frac{1}{2}\vol\circ \tau_2\;.
\end{equation}
In order to compute the $Y$--optimal hexahedra, we can first argue, that the faces are planar, then show that they are orthogonal to each other and lastly that the face areas are all equal. We leave the details to the interested reader.

\subsection{The icosahedron and the dodecahedron}\label{SecIcoDode}

Due to its triangular faces the $i$--th entry of the gradient $X$ of the mean volume is given by the normal vector to the surfaces enclosed by the link of $x_i$. However, the singularities of the gradient flow are not only regular icosahedra. In fact, the regular icosahedra are not maximal with respect to the volume function. We can vary the gradient flow by adding the normal vector to the surfaces enclosed by the link of the point opposite to $x_i$. More precisely, the resulting gradient vector $Y$ is given by Equation \eqref{EqGradIco}, where only the first entry is shown and the other entries can be found by rotation symmetry. Experiments have shown that the $Y$--optimal singularities are regular icosahedra and regular great icosahedra.\begin{equation}\label{EqGradIco}
Y_x = \frac{1}{6}\begin{pmatrix}
\nu	(2,\ldots,6)\\
\quad + \nu (7,\ldots,11)\\
\quad + \nu(5,7,3,10,6,8,4,11,2,9)\\
\vdots
\end{pmatrix}
 \quad \text{for \quad
\raisebox{-1cm}{\def\svgwidth{3.5cm}
%% Creator: Inkscape inkscape 0.48.3.1, www.inkscape.org
%% Accompanies image file 'regico.pdf' (pdf, eps, ps)
\begingroup%
  \makeatletter%
  \providecommand\color[2][]{%
    \errmessage{(Inkscape) Color is used for the text in Inkscape, but the package 'color.sty' is not loaded}%
    \renewcommand\color[2][]{}%
  }%
  \providecommand\transparent[1]{%
    \errmessage{(Inkscape) Transparency is used (non-zero) for the text in Inkscape, but the package 'transparent.sty' is not loaded}%
    \renewcommand\transparent[1]{}%
  }%
  \providecommand\rotatebox[2]{#2}%
  \ifx\svgwidth\undefined%
    \setlength{\unitlength}{296.59178467bp}%
    \ifx\svgscale\undefined%
      \relax%
    \else%
      \setlength{\unitlength}{\unitlength * \real{\svgscale}}%
    \fi%
  \else%
    \setlength{\unitlength}{\svgwidth}%
  \fi%
  \global\let\svgwidth\undefined%
  \global\let\svgscale\undefined%
  \makeatother%
  \begin{picture}(1,0.61876422)%
    \put(0,0){\includegraphics[width=\unitlength]{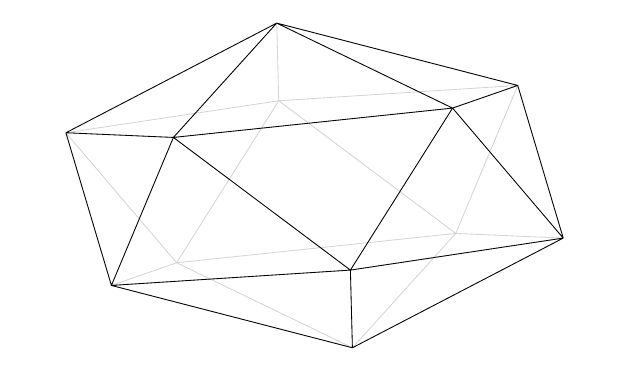}}%
    \put(0.42914381,0.59417234){\color[rgb]{0,0,0}\makebox(0,0)[lb]{\smash{$x_1$}}}%
    \put(-0.00268677,0.38583511){\color[rgb]{0,0,0}\makebox(0,0)[lb]{\smash{$x_2$}}}%
    \put(0.15587654,0.34922872){\color[rgb]{0,0,0}\makebox(0,0)[lb]{\smash{$x_3$}}}%
    \put(0.59245545,0.37986247){\color[rgb]{0,0,0}\makebox(0,0)[lb]{\smash{$x_4$}}}%
    \put(0.85120451,0.47715833){\color[rgb]{0,0,0}\makebox(0,0)[lb]{\smash{$x_5$}}}%
    \put(0.38707457,0.44845117){\color[rgb]{0,0,0}\makebox(0,0)[lb]{\smash{$\cg x_6$}}}%
    \put(0.09037042,0.115526){\color[rgb]{0,0,0}\makebox(0,0)[lb]{\smash{$x_7$}}}%
    \put(0.57781286,0.13537059){\color[rgb]{0,0,0}\makebox(0,0)[lb]{\smash{$x_8$}}}%
    \put(0.9055362,0.19259212){\color[rgb]{0,0,0}\makebox(0,0)[lb]{\smash{$x_9$}}}%
    \put(0.70304516,0.24962087){\color[rgb]{0,0,0}\makebox(0,0)[lb]{\smash{$\cg x_{10}$}}}%
    \put(0.24353912,0.20473001){\color[rgb]{0,0,0}\makebox(0,0)[lb]{\smash{$\cg x_{11}$}}}%
    \put(0.51134345,0.0076336){\color[rgb]{0,0,0}\makebox(0,0)[lb]{\smash{$x_{12}$}}}%
  \end{picture}%
\endgroup%
}.}
\end{equation}

If $X$ is the gradient of the mean volume for the dodecahedra, we have observed experimentally, that the $X$--optimal singularities are regular dodecahedra.

\section{Summary and Outlook}\label{SecOutlook}

We have shown that the mean volume function or a variation of it can be viewed as a new, simple and geometrically intuitive regularity and quality criterion for volume elements, whose gradient flow improves the quality measured by this criterion. Numerical tests and partial proofs show that the flow of the gradient vector field and its variations efficiently optimize and untangle the four polyhedron types, which are relevant for the finite element method. In the case of the tetrahedron, the volume function is a Morse-Bott function on $S^8$, whose gradient yields the efficient smoothing algorithm introduced in \cite{VartziotisWipperSchwald2009}, with one maximal and one minimal 2--sphere, as well as another set of index $(2,2)$ singularities at level 0 homeomorphic to $S^2 \times S^2$. The mean ratio quality measure 
has successfully been used in order to show the effectiveness of the GETMe smoothing. The results of this paper suggest, however, that the effectiveness of GETMe would be better measured with respect to the mean volume quality measure.

\begin{figure}[ht]
\centering
\def\svgwidth{3.4cm}
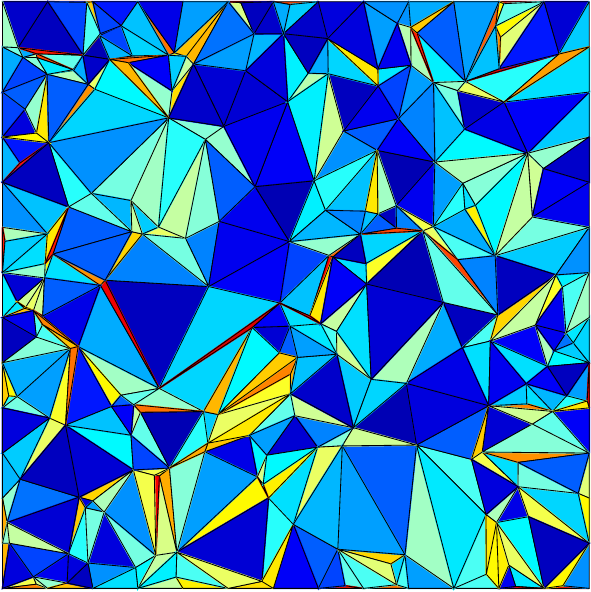
\def\svgwidth{3.4cm}
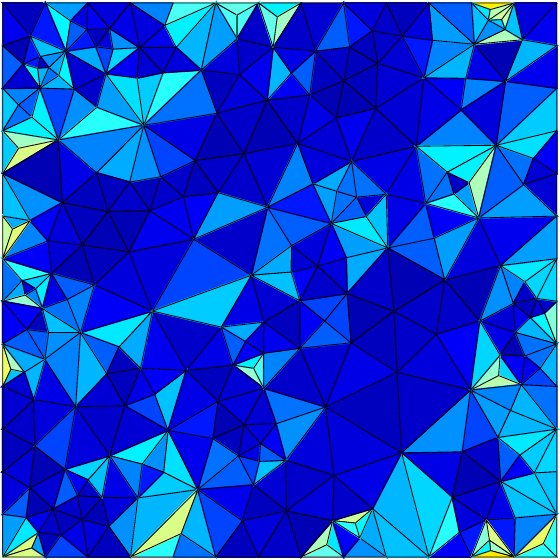
\def\svgwidth{3.4cm}
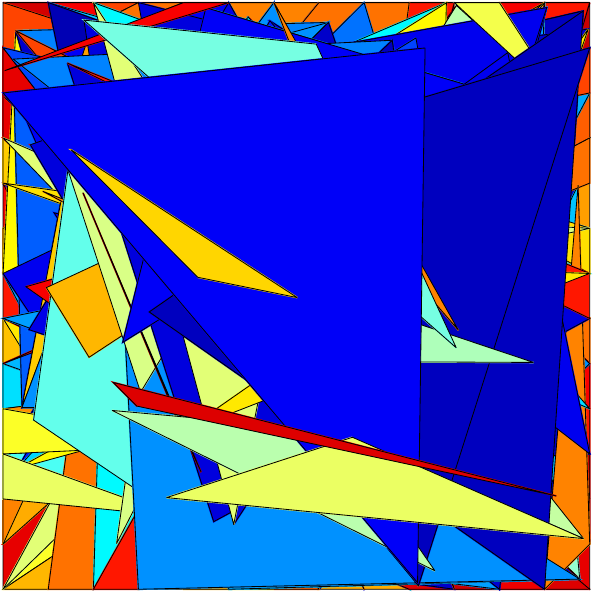
\def\svgwidth{3.4cm}
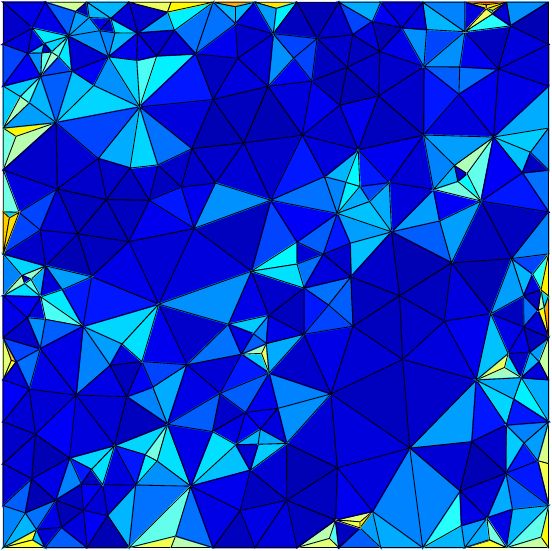\\[1.5ex]
\def\svgwidth{11.4cm}
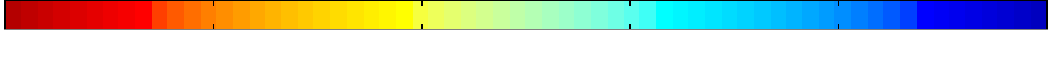
\caption{From left to right: The initial mesh with a mesh optimizing $iq(E)$, a random mesh with the same connectivity as the left mesh with a mesh optimizing $q(E)$. Each element is colored according to its mean ratio quality number.\label{fig_square_exple}}
\end{figure}

In \cite{VartziotisHimpel2013b} we have presented an idea for constructing global optimization-based GETMe smoothing algorithms. We will follow up on this idea in a future publication by considering other global quality measures. For example, if $E$ is the set of volume elements for a mesh and $x_e$ are the coordinates of the volume element $e \in E$, then we define\begin{align*}%\label{eq:qualitymeasures}
\iq(E) &= \prod_{e\in E} \iq(x_e),\quad \text{where } \iq(x) = C \frac{\vol(x)}{\area(x)^{3/2}},\\
\tag*{and} q(E) &= \sum_{e\in E} q(x_e), \quad\text{where } q(x) = \vol(x) - \frac{1}{C}\area(x)^{3/2}.
\end{align*}
with $C>0$ satisfying $\max_x\iq(x) =1$. It turns out that $q$ is concave on inner nodes and can be used for convex optimization, while giving efficient GETMe smoothing procedures. A preliminary test of this idea on a planar triangle mesh using the analogous quality function corresponding to the two-dimensional isoperimetric quotient shows that this is a promising approach. As we can see in Figure~\ref{fig_square_exple}, while the results from optimizing $\iq(E)$ and $q(E)$ look similar, the latter makes smaller triangles worse and bigger triangles better.

\

\noindent{\bf Acknowledgements.} We thank Joachim Wipper from TWT GmbH Science \& Innovation, Department for Mathematical Research \& Services, for fruitful discussions and some numerical tests.

\bibliographystyle{siam}
\bibliography{references.bib}

\end{document}